\documentclass{amsart}
\usepackage{mathrsfs}
\usepackage{mathrsfs}
\usepackage{amsfonts}
\usepackage{cases}
\usepackage{latexsym}
\usepackage{amsmath}
\usepackage[all]{xy}
\usepackage{stmaryrd}
\usepackage{amsfonts}
\usepackage{amsmath,amssymb,amscd,bbm,amsthm,mathrsfs,dsfont}
\usepackage[notref, notcite]{showkeys}

\usepackage[numbers,sort&compress]{natbib}
\usepackage[pagebackref]{hyperref}

\usepackage{color, hyperref}
\definecolor{blue}{rgb}{1.0,0.0,0.9}
\hypersetup{colorlinks, breaklinks,
            linkcolor=blue,urlcolor=blue,
            anchorcolor=blue, citecolor=blue}

\usepackage{fancyhdr}
\usepackage{amsxtra,ifthen}
\usepackage{verbatim}

\numberwithin{equation}{section}

\usepackage{fancyhdr}
\usepackage{amsxtra,ifthen}
\usepackage{verbatim}

\usepackage[numbers,sort&compress]{natbib}
\usepackage[pagebackref]{hyperref}

\theoremstyle{plain}
\newtheorem{theorem}{Theorem}[section]
\newtheorem{proposition}[theorem]{Proposition}
\newtheorem{lemma}[theorem]{Lemma}
\newtheorem{corollary}[theorem]{Corollary}

\newtheorem{question}[theorem]{Question}

\theoremstyle{definition}
\newtheorem{definition}[theorem]{Definition}

\begin{document}

\title[$k$-Commuting Mappings of Generalized Matrix Algebras]
{$k$-Commuting Mappings of Generalized Matrix Algebras}

\author{Yanbo Li, Feng Wei and Ajda Fo\v sner}

\address{Li: School of Mathematics and Statistics, Northeastern University at Qinhuangdao, Qinhuangdao, 066004, P. R. China.}

\email{\href{liyanbo707@163.com}{liyanbo707@163.com}}

\address{Wei: School of Mathematics and Statistics, Beijing Institute of
Technology, Beijing, 100081, P. R. China}

\email{\href{daoshuo@hotmail.com}{daoshuo@hotmail.com}}\email{\href{daoshuowei@gmail.com}{daoshuowei@gmail.com}}

\address{Fo\v sner: Faculty of Management, University of Primorska,
Cankarjeva 5, SI-6104 Koper, Slovenia}

\email{\href{ajda.fosner@fm-kp.si}{ajda.fosner@fm-kp.si}}

\begin{abstract}
In this paper we will study $k$-commuting mappings of generalized
matrix algebras. The general form of arbitrary $k$-commuting mapping
of a generalized matrix algebra is determined. It is shown that
under mild assumptions, every $k$-commuting mapping of a
generalized matrix algebra takes a certain form which is said to be
proper. A number of applications related to $k$-commuting mappings
are presented. These results extend the existing works of Cheung, Du
and Wang \cite{Cheung2, DuWang1} to the case of generalized matrix
algebras.
\end{abstract}

\date{\today}

\subjclass[2000]{16R60, 16W10, 15A78}

\keywords{$k$-Commuting mapping, generalized matrix algebra}

\thanks{ This work is partially supported by the Training Program of International Exchange and
Cooperation of the Beijing Institute of Technology. The work of the first author is supported by the Natural Science Foundation of Hebei
Province, China (A2017501003) and the Science and Technology support
program of Northeastern University at Qinhuangdao (No. XNK201601).}

\maketitle

\section{Introduction}
\label{xxsec1}

Let $\mathcal{R}$ be a commutative ring with identity element and
$\mathcal{A}$ a unital associative $\mathcal{R}$-algebra. For arbitrary
elements $a, b\in \mathcal{A}$, we set $[a, b]_0=a$, $[a, b]_1=ab-ba$, and
inductively $[a, b]_k = [[a, b]_{k-1}, b]$, where $k$ is a fixed
positive integer. Denote by $\mathcal{Z(A)}$ the center of $\mathcal{A}$
and define
$$
{\mathcal Z(A)}_k =\{\hspace{2pt}a\in \mathcal{A}\hspace{2pt}|\hspace{2pt} [a,
x]_k=0, \hspace{2pt} \forall x\in \mathcal{A}\hspace{2pt}\}.
$$
Clearly, $\mathcal{Z(A)}_1=\mathcal{Z(A)}$. An
$\mathcal{R}$-linear mapping $\Theta\colon \mathcal{A}\longrightarrow \mathcal{A}$ is said
to be $k$-\textit{commuting} on $\mathcal{A}$ if $[\Theta(a), a]_k = 0$ for
all $a\in \mathcal{A}$. In particular, an $\mathcal{R}$-linear mapping
$\Theta\colon \mathcal{A}\longrightarrow \mathcal{A}$ is called \textit{commuting} on $\mathcal{A}$ if
$[\Theta(a), a]=0$ for all $a\in \mathcal{A}$. When we investigate a
$k$-commuting mapping $\Theta$ of an algebra $\mathcal{A}$, the principal task
is to describe its form. Let $\Theta$ be a $k$-commuting mapping of
an $\mathcal{R}$-algebra $\mathcal{A}$. Then $\Theta$ will be called
\textit{proper} if it has the form
$$
\Theta(a)=\lambda a+\zeta(a)\eqno(\clubsuit)
$$
for all $a\in \mathcal{A}$, where $\lambda\in \mathcal{Z(A)}$ and $\zeta\colon \mathcal{A}\longrightarrow \mathcal{Z(A)}$ is an $\mathcal{R}$-linear
mapping. The concept of commuting mappings is closely related to that of biderivations. Recall that
a bilinear mapping $\Upsilon\colon \mathcal{A}\times \mathcal{A}\longrightarrow \mathcal{A}$ is
called a \textit{biderivation} if it is a derivation with respect to both components,
meaning that $\Upsilon(ab, c)=a\Upsilon(b, c)+\Upsilon(a, c)b$ and $\Upsilon(a, bc)=b\Upsilon(a, c)+\Upsilon(a, b)c$ for all $a, b, c\in \mathcal{A}$.
If $\mathcal{A}$ is a noncommutative algebra, then the mapping
$\Upsilon(a, b)=\lambda[a, b]$ (where $\lambda\in \mathcal{Z(A)}, \forall a, b\in \mathcal{A}$) is a basic
example of biderivation, which is usually called
\textit{inner biderivations}. Every commuting mapping $\Theta\colon \mathcal{A}\longrightarrow \mathcal{A}$
gives rise to a biderivation of $\mathcal{A}$. Namely, linearizing the relation $[\Theta(a), a]=0$ yields that
$[\Theta(a), b]=[a, \Theta(b)]$ for all $a, b\in \mathcal{A}$. Thus the resulting mapping
$(a, b)\longmapsto [\Theta(a), b]$ is a biderivation.

Commuting mappings and biderivations are currently active branches in the theory of additive mappings of
noncommutative algebras. Bre\v{s}ar \cite{Bresar2} showed that an $\mathcal{R}$-linear mapping
$\Theta$ of a prime algebra $\mathcal{A}$ is commuting if and only if it has
the form $(\clubsuit)$. This result gives rise to the study of
various more general problems, and eventually comes into the theory
of functional identities \cite{BresarChebotarMartindale}. We also
encourage the reader to read the elegant survey paper
\cite{Bresar4}, in which the author presented a full and detailed
account for the theory of commuting mappings. Furthermore,
Bre\v{s}ar considered linear mappings with Engel condition on prime
algebras, and especially studied $2$-commuting and $k$-commuting
mappings of prime algebras in \cite{Bresar1, Bresar2}. He observed
that $2$-commuting and $k$-commuting mappings of certain prime
algebras are commuting. Zhang et al. \cite{ZhangFengLiWu} showed that every
linear biderivation of nest algebras on a complex separable Hilbert space $\textbf{H}$ is an
inner biderivation if and only if ${\rm dim}\ 0_+\neq 1$ or ${\rm dim}\ \textbf{H}_{-}^{\perp}\neq 1$.
Zhao et al. \cite{ZhaoWangYao} used the results in \cite{ZhangFengLiWu} to prove that every
biderivation of an upper triangular matrix algebra is a sum of an inner biderivation and a special
biderivation, which they call an extremal biderivation. Benkovi\v{c} \cite{Benkovic1} defines the concept of an extremal biderivation,
and proves that under certain conditions a biderivation of a triangular
algebra $\mathcal{A}$ is a sum of an extremal and an inner biderivation.
The same statement still holds for the case of generalized matrix algebras \cite{DuWang3}.

Cheung initially started to study commuting mappings of matrix
algebras in his powerful works \cite{Cheung1, Cheung2}. He determined the
class of triangular algebras for which every commuting
mapping is proper. In \cite{BenkovicEremita} Benkovi\v{c} and Eremita studied commuting
traces of bilinear mappings on triangular algebras. They gave mild
conditions under which arbitrary commuting trace of a triangular
algebra is proper. The authors applied the obtained results to the
study of Lie isomorphisms and that of commutativity preserving
mappings. Du and Wang \cite{DuWang1} proved that under certain
conditions, each $k$-commuting mapping on a triangular algebra $\mathcal{A}$ is
proper. More recently, Li, Liang, Wei and Xiao \cite{FosnerLiangWei, LivanWykWei, LiWei, LiXiao,
LiangWeiXiaoFosner, XiaoWei1, XiaoWei2, XiaoWeiFosner} jointly investigated linear mappings of generalized
matrix algebras, such as derivations, Jordan derivations, Lie derivations, commuting
mappings and semi-centralizing mappings. Our main purpose is to
develop the theory of linear mappings of triangular algebras to the
case of generalized matrix algebras, which has a much broader
background. In \cite{XiaoWei1}, Xiao and Wei extended the main
results of \cite{Cheung2} to the case of generalized matrix algebras.
They described the general form of arbitrary commuting mapping of a
generalized matrix algebra and provided several sufficient
conditions which enable the commuting mappings to be proper. Li and
Wei \cite{LiWei, LivanWykWei} considered semi-centralizing mappings
of generalized matrix algebras and many ring-theoretic aspect
results were extended to the case of generalized matrix algebras via
complicated matrix computations. Franca \cite{Franca1, Franca2, Franca3, Franca4,
Franca5, Franca6, Franca7, FrancaLouza7} considered commuting mappings
on certain subset of the full matrix algebra $M_n(\mathbb{K})$ over an arbitrary filed $\mathbb{K}$.
The involved subsets include the subset of all invertible matrices, the subset of all
singular matrices and the subset of all rank-$k$ matrices. He observed that
every commuting mapping on these subsets usually has the so-called proper form.
Motivating by Franca's work, Liu \cite{Liu} characterized centralizing mappings on the
above-mentioned subsets and got some analogous results. Xu and Yi gave the
forms of commuting mappings of the aforementioned subsets got some analogous results \cite{XuYi}.

This paper is devoted to the study of $k$-commuting mappings of
generalized matrix algebras. We will describe the general form of
arbitrary $k$-commuting mapping of a $2$-torsion free generalized
matrix algebra and provide a sufficient condition which enables each
commuting mapping to be proper. Our work extends the main results of
\cite{Cheung2, DuWang1} to the case of generalized matrix algebras
and also give the corresponding $k$-commuting version of
\cite[Theorem 3.6]{XiaoWei1}.

\section{Generalized Matrix Algebras and Examples}\label{xxsec2}

Let us begin with the definition of generalized matrix algebras
given by a Morita context. Let $\mathcal{R}$ be a commutative ring
with identity. A \textit{Morita context} consists of two unital
$\mathcal{R}$-algebras $A$ and $B$, two bimodules $_AM_B$ and
$_BN_A$, and two bimodule homomorphisms called the pairings
$\Phi_{MN}: M\underset {B}{\otimes} N\longrightarrow A$ and
$\Psi_{NM}: N\underset {A}{\otimes} M\longrightarrow B$ satisfying
the following commutative diagrams:
$$
\xymatrix{ M \underset {B}{\otimes} N \underset{A}{\otimes} M
\ar[rr]^{\hspace{8pt}\Phi_{MN} \otimes I_M} \ar[dd]^{I_M \otimes
\Psi_{NM}} && A
\underset{A}{\otimes} M \ar[dd]^{\cong} \\  &&\\
M \underset{B}{\otimes} B \ar[rr]^{\hspace{10pt}\cong} && M }
\hspace{6pt}{\rm and}\hspace{6pt} \xymatrix{ N \underset
{A}{\otimes} M \underset{B}{\otimes} N
\ar[rr]^{\hspace{8pt}\Psi_{NM}\otimes I_N} \ar[dd]^{I_N\otimes
\Phi_{MN}} && B
\underset{B}{\otimes} N \ar[dd]^{\cong}\\  &&\\
N \underset{A}{\otimes} A \ar[rr]^{\hspace{10pt}\cong} && N
\hspace{2pt}.}
$$
Let us write this Morita context as $(A, B, M, N, \Phi_{MN},
\Psi_{NM})$. We refer the reader to \cite{Morita} for the basic
properties of Morita contexts. If $(A, B, M, N,$ $ \Phi_{MN},
\Psi_{NM})$ is a Morita context, then the set
$$
\left[
\begin{array}
[c]{cc}%
A & M\\
N & B\\
\end{array}
\right]=\left\{ \left[
\begin{array}
[c]{cc}%
a& m\\
n & b\\
\end{array}
\right] \vline a\in A, m\in M, n\in N, b\in B \right\}
$$
is an $\mathcal{R}$-algebra under matrix-like addition and
matrix-like multiplication, where at least one of
$M$ and $N$ is non-zero. Such an $\mathcal{R}$-algebra is
usually called a \textit{generalized matrix algebra} of order $2$
and is denoted by
$$
\mathcal{G}=\mathcal{G}(A, M, N, B)=\left[
\begin{array}
[c]{cc}%
A & M\\
N & B\\
\end{array}
\right].
$$
In a similar way, one can define a generalized matrix algebra of
order $n>2$. It was shown that up to isomorphism, arbitrary
generalized matrix algebra of order $n$ $(n\geq 2)$ is a generalized
matrix algebra of order 2 \cite[Example 2.2]{LiWei}. If one of the
modules $M$ and $N$ is zero, then $\mathcal{G}$ exactly degenerates
to an \textit{upper triangular algebra} or a \textit{lower
triangular algebra}. In this case, we denote the resulted upper
triangular algebra (resp. lower triangular algebra) by
$$\mathcal{T^U}=\mathcal{T}(A, M, B)=
\left[
\begin{array}
[c]{cc}%
A & M\\
O & B\\
\end{array}
\right]   \hspace{8pt} \left({\rm resp.} \hspace{4pt} \mathcal{T_L}=\mathcal{T}(A, N, B)=
\left[
\begin{array}
[c]{cc}%
A & O\\
N & B\\
\end{array}
\right]\right)
$$
Note that our current generalized matrix algebras contain those
generalized matrix algebras in the sense of Brown \cite{Brown} as
special cases. Let $\mathcal{M}_n(\mathcal{R})$ be the full matrix
algebra consisting of all $n\times n$ matrices over $\mathcal{R}$.
It is worth to point out that the notion of generalized matrix
algebras efficiently incorporates triangular algebras and full matrix
algebras together. A distinctive feature of our systematic work
is to deal with all questions related to (non-)linear mappings of
triangular algebras and of full matrix algebras under a rigorous unified framework, which is the admired generalized matrix algebras frame, see
\cite{LivanWykWei, LiWei, LiXiao, LiangWeiXiaoFosner, XiaoWei1, XiaoWei2}.

Let us list some classical examples of generalized matrix algebras
which will be revisited in the sequel (Section \ref{xxsec4} and
Section \ref{xxsec5}). Since these examples are ubiquitous in  noncommutative background and operator theory, we just state their title without any details.
\begin{enumerate}
\item[(1)] Unital algebras with nontrivial
idempotents, such as (semi-)prime algebras with nontrivial idempotents;
\item[(2)] Full matrix algebras;
\item[(3)] Inflated algebras;
\item[(4)] Triangular algebras;
\item[(5)] Quasitilted algebras;
\item[(6)] von Neumann algebra on Hilbert spaces;
\item[(7)] Nest algebras on Hilbert spaces;
\item[(8)] Standard operator algebras on Banach spaces.
\end{enumerate}

\noindent{Special} attention is paid to the unital algebras with nontrivial idempotents. The following result demonstrates the equivalence between
the class of generalized matrix algebras and the family of unital algebras with nontrivial idempotents.

\begin{proposition}{\rm \cite[Proposition 2.1]{XiaoWei1}}\label{xxsec2.1}
A unital algebra $\mathscr{A}$ is a generalized matrix algebra if and only if there exists an
idempotent $e\in \mathscr{A}$ such that $e\mathscr{A}(1-e)\neq 0$.
\end{proposition}

\noindent These generalized matrix algebras regularly appear in the theory of associative algebras and
noncommutative Noetherian algebras in the most diverse situations, which is due to its powerful persuasiveness
and intuitive illustration effect. However, people pay less attention to the linear mappings of generalized matrix
algebras. It was Krylov who initiated the study of linear mappings on generalized matrix algebras from the classifying point
of view \cite{Krylov1}. Since then many articles are devoted to this topic, and a number of interesting
results are obtained (see \cite{AnhWyk1, AnhWyk2, AnhBirkenmeierWyk, Benkovic1, BenkovicGrasic, BenkovicSirovnik,
BobocDascalescuWyk, DuWang2, DuWang3, LivanWykWei, LiWei, LiXiao, LiangWeiXiaoFosner, WangWang, XiaoWei1, XiaoWei2, XiaoWeiFosner}).
Nevertheless, it leaves so much to be desired. It seems that the essential difference between triangular algebras and
generalized matrix algebras lie in the right upper ``corner" or left lower ``corner". Such a bit difference clearly increases
the complexity of construction and computation, which is reflected in the maximal ring of quotients \cite{Muller, Stenstrom}
and the modules over generalized matrix algebras \cite{KrylovTuganbaev1}. The representation theory, homological
behavior, $K$-theory of generalized matrix algebras are intensively investigated by Krylov and his coauthors in
\cite{Krylov1, Krylov2, Krylov3, KrylovTuganbaev1, KrylovTuganbaev2, KrylovTuganbaev3}. We will propose
some open questions concerning linear mappings and functional identities of generalized matrix algebras in
Section \ref{xxsec5} of this article. Therefore, generalized matrix algebras are indeed one class of great potential
and inspiring associative algebras. We can never emphasize on the importance of generalized matrix algebras too much.

\section{$k$-Commuting Mappings of Generalized Matrix Algebras}
\label{xxsec3}

Throughout this section, we denote the generalized matrix algebra
of order $2$ originating from a Morita context $(A, B, M, N,
\Phi_{MN}, \Psi_{NM})$ by
$$
\mathcal{G}:=\left[
\begin{array}
[c]{cc}%
A & M\\
N & B
\end{array}
\right] ,
$$
where at least one of the two bimodules $M$ and $N$ is distinct from
zero. We always assume that $M$ is faithful as a left $A$-module and
also as a right $B$-module, but no any constraint conditions on $N$.
By \cite[Section 2.2]{LiWei} we know that every generalized matrix
algebra of order $n (n>2)$ is isomorphic to a generalized matrix
algebras of order $2$. In view of this fact and technical
considerations, only generalized matrix algebras of order $2$ are
considered in this section. Let $k$ be a fixed positive integer with
$k\geq 2$. An $\mathcal{R}$-linear mapping $\Theta:
\mathcal{G}\longrightarrow \mathcal{G}$ is called
$k$-\textit{commuting} if
$$
\left[\Theta\left(\left[
\begin{array}
[c]{cc}%
a & m\\
n & b
\end{array}
\right]\right), \left[
\begin{array}
[c]{cc}%
a & m\\
n & b
\end{array}
\right]\right]_k=\left[
\begin{array}
[c]{cc}%
0 & 0\\
0 & 0
\end{array}
\right]
$$
for all $\left[
\begin{array}
[c]{cc}%
a & m\\
n & b
\end{array}
\right]\in \mathcal{G}$.

The center of $\mathcal{G}$ is
$$
\mathcal{Z(G)}=\left\{ \left[
\begin{array}
[c]{cc}%
a & 0\\
0 & b
\end{array}
\right] \vline \hspace{3pt} am=mb, \hspace{3pt} na=bn,\ \forall\
m\in M, \hspace{3pt} \forall n\in N \right\}.
$$
Indeed, by \cite[Lemma 1]{Krylov1} it follows that the center
$\mathcal{Z(G)}$ consists of all diagonal matrices $
\left[\smallmatrix a & 0\\
0 & b
\endsmallmatrix \right]$,
where $a\in \mathcal{Z}(A)$, $b\in \mathcal{Z}(B)$ and $am=mb$,
$na=bn$ for all $m\in M, n\in N$. However, in our situation where
$M$ is faithful as a left $A$-module and also as a right $B$-module,
the conditions that $a\in \mathcal{Z}(A)$ and $b\in \mathcal{Z}(B)$
become redundant and can be deleted. Indeed, if $am=mb$ for all
$m\in M$, then for arbitrary element $a^\prime \in A$ we get
$$
(aa^\prime-a^\prime a)m=a(a'm)-a'(am)=(a'm)b-a'(mb)=0.
$$
The assumption that $M$ is faithful as a left $\mathcal{A}$-module
leads to $aa'-a'a=0$ and hence $a\in \mathcal{Z}(A)$. Likewise, we
also have $b\in \mathcal{Z}(B)$.

Let us define two natural $\mathcal{R}$-linear projections
$\pi_A:\mathcal{G}\rightarrow A$ and $\pi_B:\mathcal{G}\rightarrow
B$ by
$$
\pi_A: \left[
\begin{array}
[c]{cc}%
a & m\\
n & b\\
\end{array}
\right] \longmapsto a \quad \text{and} \quad \pi_B: \left[
\begin{array}
[c]{cc}%
a & m\\
n & b\\
\end{array}
\right] \longmapsto b.
$$
By the above paragraph, it is not difficult to see that $\pi_A
\left(\mathcal{Z(G)}\right)$ is a subalgebra of $\mathcal{Z}(A)$ and
that $\pi_B\left(\mathcal{Z(G)}\right)$ is a subalgebra of
$\mathcal{Z}(B)$.
Given an element $a\in\pi_A(\mathcal{Z(G)})$, if $\left[\smallmatrix a & 0\\
0 & b
\endsmallmatrix \right], \left[\smallmatrix a & 0\\
0 & b^\prime
\endsmallmatrix \right] \in \mathcal{Z(G)}$, then we have $am=mb=mb'$ for
all $m\in M$. Since $M$ is faithful as a right $B$-module,
$b=b^\prime$. That implies there exists a unique
$b\in\pi_B(\mathcal{Z(G)})$, which is denoted by $\varphi(a)$, such
that $
\left[\smallmatrix a & 0\\
0 & b
\endsmallmatrix \right] \in \mathcal{Z(G)}$. It is easy to
verify that the map $\varphi:\pi_A(\mathcal{Z(G)})\longrightarrow
\pi_B(\mathcal{Z(G)})$ is an algebraic isomorphism such that
$am=m\varphi(a)$ and $na=\varphi(a)n$ for all $a\in
\pi_A(\mathcal{Z(G)}), m\in M, n\in N$.

The following result is a natural extension of \cite[Lemma
2.1]{DuWang1}, which is indispensable for the proof of our main
result.

\begin{lemma}\label{xxsec3.1}
Let $n$ be a positive integer and $A$ be a unital associative
ring. For a left $A$-module $M$, if $\alpha: A\rightarrow M$ is a
mapping such that $\alpha(a + 1) = \alpha(a)$ and $a^n\alpha(a) =
0$ for all $a\in A$, then $\alpha = 0$. Similarly, for a right
$A$-module $M^\prime$, a mapping $\beta: A\rightarrow M^\prime$ is
zero if $\beta(a + 1) = \beta(a)$ and $\beta(a)a^n = 0$ for all
$a\in A$.
\end{lemma}

Before proving our main theorem, we describe the general form of
arbitrary $k$-commuting mapping on the generalized matrix algebra
$\mathcal{G}$.

\begin{proposition}\label{xxsec3.2}
Let $\Theta$ be a $k$-commuting mapping of $\mathcal{G}$. Then
$\Theta$ is of the form
$$
\begin{aligned}
& \Theta\left(\left[
\begin{array}
[c]{cc}%
a & m\\
n & b\\
\end{array}
\right]\right)\\
=& \left[
\begin{array}
[c]{cc}%
\delta_1(a)+\delta_2(m)+\delta_3(n)+\delta_4(b) & \tau_2(m)\\
\nu_3(n) & \mu_1(a)+\mu_2(m)+\mu_3(n)+\mu_4(b)\\
\end{array}
\right] ,
\end{aligned}  \eqno(\bigstar)
$$
where
$$
\begin{aligned} \delta_1:& A \longrightarrow A, &
\delta_2: & M\longrightarrow {\mathcal Z}(A)_k, & \delta_3: & N
\longrightarrow
{\mathcal Z}(A)_k, & \delta_4: & B\longrightarrow {\mathcal Z}(A)_k,\\
\mu_1: & A\longrightarrow {\mathcal Z}(B)_k, & \mu_2: &
M\longrightarrow {\mathcal Z}(B)_k, & \mu_3: & N\longrightarrow
{\mathcal Z}(B)_k, & \mu_4:
&B\longrightarrow B\\
& & \tau_2:&  M\longrightarrow M, & \nu_3: & N\longrightarrow N &&&
\end{aligned}
$$
are all $\mathcal{R}$-linear mappings satisfying the following
conditions:
\begin{enumerate}
\item[{\rm(1)}] $\delta_1$ is a $k$-commuting mapping of $A$
and $\delta_1(1)\in {\mathcal Z}(A)_k;$
\item[{\rm(2)}] $\mu_4$ is a $k$-commuting mapping of $B$ and
$\mu_4(1)\in {\mathcal Z}(B)_k;$
\item[{\rm(3)}]
$(\delta_1(1)+\delta_4(1)+2\delta_2(m))m=m(\mu_1(1)+\mu_4(1)+2\mu_2(m));$
\item[{\rm(4)}]
$n(\delta_1(1)+\delta_4(1)+2\delta_3(n))=(\mu_1(1)+\mu_4(1)+2\mu_3(n))n;$
\item[{\rm(5)}]
$2\tau_2(m)=(\delta_1(1)-\delta_4(1))m-m(\mu_1(1)-\mu_4(1));$
\item[{\rm(6)}]
$2\nu_3(n)=n(\delta_1(1)-\delta_4(1))-(\mu_1(1)-\mu_4(1))n.$
\end{enumerate}
\end{proposition}

\begin{proof}
Suppose that the $k$-commuting mapping $\Theta$ is of the form
$$
\begin{aligned}
& \Theta\left(\left[
\begin{array}
[c]{cc}%
a & m \\
n & b \\
\end{array}
\right]\right) \\
= &\left[
\begin{array}
[c]{cc}%
\delta_1(a)+\delta_2(m)+\delta_3(n)+\delta_4(b) & \tau_1(a)+\tau_2(m)+\tau_3(n)+\tau_4(b) \\
\nu_1(a)+\nu_2(m)+\nu_3(n)+\nu_4(b) & \mu_1(a)+\mu_2(m)+\mu_3(n)+\mu_4(b) \\
\end{array}
\right]
\end{aligned} \eqno(3.1)
$$
for all $\left[\smallmatrix a & m\\
n & b
\endsmallmatrix \right]\in \mathcal{G}$, where $\delta_1,\delta_2,\delta_3,\delta_4$ are
$\mathcal{R}$-linear mappings from $A, M, N, B$ to $A$,
respectively; $\tau_1,\tau_2$, $\tau_3,\tau_4$ are
$\mathcal{R}$-linear mappings from $A, M, N, B$ to $M$,
respectively; $\nu_1,\nu_2,\nu_3,\nu_4$ are $\mathcal{R}$-linear
mappings from $A, M, N, B$ to $N$, respectively;
$\mu_1,\mu_2,\mu_3,\mu_4$ are $\mathcal{R}$-linear mappings from $A,
M, N, B$ to $B$, respectively.

For any $G\in\mathcal {G}$, we will intensively employ the equation
$$
[\Theta(G), G]_k=\left[
\begin{array}
[c]{cc}%
0 & 0 \\
0 & 0 \\
\end{array}
\right].\eqno(3.2)
$$
Taking $G=\left[\smallmatrix 1 & 0\\
0 & 0
\endsmallmatrix \right]$ into (3.1) leads to
$$\Theta(G)=\left[
\begin{array}
[c]{cc}%
\delta_1(1) & \tau_1(1) \\
\nu_1(1) & \mu_1(1) \\
\end{array}
\right]. \eqno(3.3)
$$
Combining (3.2) with (3.3) and a direct computation yields
$$
\left[
\begin{array}
[c]{cc}%
0 & 0 \\
0 & 0 \\
\end{array}
\right]=[\Theta(G),G]_k =\left[
\begin{array}
[c]{cc}%
0 & (-1)^k\tau_1(1) \\
\nu_1(1) & 0 \\
\end{array}
\right].
$$
This implies that
$$
\tau_1(1)=0, \quad \nu_1(1)=0. \eqno(3.4)
$$
Likewise, we also have
$$
\tau_4(1)=0, \quad \nu_4(1)=0 \eqno(3.5)
$$
by putting $G=\left[\smallmatrix 0 & 0\\
0 & 1
\endsmallmatrix \right]$ in (3.2).

Let us take $G=\left[\smallmatrix a & 0\\
0 & 0
\endsmallmatrix \right]$ into (3.2).
An inductive approach gives
$$
\left[
\begin{array}
[c]{cc}%
0 & 0 \\
0 & 0 \\
\end{array}
\right]=[\Theta(G),G]_k =\left[
\begin{array}
[c]{cc}%
[\delta_1(a), a]_k & (-1)^ka^k\tau_1(a) \\
\nu_1(a)a^k & 0 \\
\end{array}
\right].
\eqno(3.6)
$$
This shows that
$$
[\delta_1(a), a]_k=0 \eqno(3.7)
$$
for all $a\in A$. That is, $\delta_1$ is a $k$-commuting mapping of
$A$. Substituting $a+1$ for $a$ in (3.7) we get $[\delta_1(1),
a]_k=0$ for all $a\in A$. Therefore $\delta_1(1)\in {\mathcal
Z}(A)_k$. By (3.6) we know that $a^k\tau_1(a)=0$. In view of (3.4)
we obtain $\tau_1(a)=\tau_1(a+1)$. By Lemma \ref{xxsec3.1} it
follows that
$$
\tau_1(a)=0 \eqno(3.8)
$$
for all $a\in A$. Revisiting the relations (3.6) and (3.4) and
applying Lemma \ref{xxsec3.1} again we have
$$
\nu_1(a)=0 \eqno(3.9)
$$
for all $a\in A$.

Let us choose $G=\left[\smallmatrix 0 & 0\\
0 & b\endsmallmatrix \right]$ in (3.2). Then
$$
\left[
\begin{array}
[c]{cc}%
0 & 0 \\
0 & 0 \\
\end{array}
\right]=[\Theta(G),G]_k =\left[
\begin{array}
[c]{cc}%
0 & \tau_4(b)b^k \\
(-1)^kb^k\nu_4(b) & [\mu_4(b), b]_k \\
\end{array}
\right].\eqno(3.10)
$$
(3.10) implies that
$$
[\mu_4(b), b]_k=0 \eqno(3.11)
$$
for all $b\in B$. That is to say that $\mu_4$ is a $k$-commuting
mapping of $B$. Replacing $b$ by $b+1$ in (3.11) leads to
$[\mu_4(1), b]_k=0$ for all $b\in B$. And hence $\mu_4(1)\in
{\mathcal Z}(B)_k$. Furthermore, it follows from (3.5), (3.10) and
Lemma \ref{xxsec3.1} that
$$
\tau_4(b)=0 \quad {\rm and}\quad \nu_4(b)=0 \eqno(3.12)
$$
for all $b\in B$.

Let us choose
$G=\left[\smallmatrix a & 0\\
0 & b
\endsmallmatrix \right]$ in (3.1). In view of (3.8), (3.9) and
(3.12) we obtain
$$
\Theta(G)=\left[
\begin{array}
[c]{cc}%
\delta_1(a)+\delta_4(b) & 0 \\
0 & \mu_1(a)+\mu_4(b) \\
\end{array}
\right]. \eqno(3.13)
$$
Taking (3.13) into (3.2) yields
$$
\begin{aligned}
\left[
\begin{array}
[c]{cc}%
0 & 0 \\
0 & 0 \\
\end{array}
\right]&=[\Theta(G),G]_k \\
&=\left[
\begin{array}
[c]{cc}%
[\delta_1(a)+\delta_4(b), a]_k & 0 \\
0 & [\mu_1(a)+\mu_4(b), b]_k \\
\end{array}
\right]\\
&=\left[
\begin{array}
[c]{cc}%
[\delta_1(a), a]_k+[\delta_4(b), a]_k & 0 \\
0 & [\mu_1(a), b]_k+[\mu_4(b), b]_k \\
\end{array}
\right].
\end{aligned}
\eqno(3.14)
$$
Note that $\delta_1$ and $\mu_4$ are $k$-commuting mappings of $A$
and $B$, respectively. Thus $[\delta_1(a), a]_k=0$ for all $a\in A$
and $[\mu_4(b), b]_k=0$ for all $b\in B$. Then (3.14) shows that
$[\delta_4(b), a]_k=0$ and $[\mu_1(a), b]_k=0$ for all $a\in A, b\in
B$. That is, $\delta_4(b)\in {\mathcal Z}(A)_k$ for all $b\in B$ and
$\mu_1(a)\in {\mathcal Z}(B)_k$ for all $a\in A$.

Let us put $G=\left[\smallmatrix 1 & m\\
0 & 0\endsmallmatrix \right]$ in (3.2) and denote by
$$
[\Theta(G),G]_i =X_i=\left[
\begin{array}
[c]{cc}%
X_i(11) & X_i(12) \\
X_i(21) & X_i(22) \\
\end{array}
\right]
$$
for each $0\leq i <k$. Then
$$
\begin{aligned}
X_{i+1}&=\left[
\begin{array}
[c]{cc}%
X_{i+1}(11) & X_{i+1}(12) \\
X_{i+1}(21) & X_{i+1}(22) \\
\end{array}
\right]\\
&=[X_i, G]\\
&=\left[\left[
\begin{array}
[c]{cc}%
X_i(11) & X_i(12) \\
X_i(21) & X_i(22) \\
\end{array}
\right],\left[
\begin{array}
[c]{cc}%
1 & m \\
0 & 0 \\
\end{array}
\right]\right]\\
&=\left[
\begin{array}
[c]{cc}%
-mX_i(11) & X_i(11)m-mX_i(22)-X_i(12) \\
X_i(21) & X_i(21)m \\
\end{array}
\right].
\end{aligned}
$$
This gives $X_{i+1}(21)=X_{i}(21)$ and hence
$$
X_k(21)=X_0(21)=\nu_2(m).\eqno(3.15)
$$
Note that the fact $X_k=0$. Then (3.15) implies that $\nu_2(m)=0$
for all $m\in M$. Therefore
$$
X_0=\left[
\begin{array}
[c]{cc}%
\delta_1(1)+\delta_2(m) & \tau_2(m) \\
0 & \mu_1(1)+\mu_2(m) \\
\end{array}
\right]
$$
and
$$
\begin{aligned}
X_1&=[X_0, G]\\
&=\left[
\begin{array}
[c]{cc}%
0 & \delta_1(1)m+\delta_2(m)m-\tau_2(m)-m\mu_1(1)-m\mu_2(m) \\
0 & 0 \\
\end{array}
\right]
\end{aligned}\eqno(3.16)
$$
Applying inductive computations we assert that for each $i>0$,
$X_i=(-1)^{i-1}X_1$ and hence $X_k=(-1)^{k-1}X_1$. This proves that
that $X_1=0$. By (3.16) we have
$$
\tau_2(m)=\delta_1(1)m+\delta_2(m)m-m\mu_1(1)-m\mu_2(m).\eqno(3.17)
$$
Likewise, we put $G=\left[\smallmatrix 0 & m\\
0 & 1\endsmallmatrix \right]$ in (3.2) and get
$$
\tau_2(m)=m\mu_4(1)+m\mu_2(m)-\delta_4(1)m-\delta_2(m)m.\eqno(3.18)
$$
Combining (3.17) with (3.18) leads to
$$
(\delta_1(1)+\delta_4(1)+2\delta_2(m))m=m(\mu_1(1)+\mu_4(1)+2\mu_2(m))
$$
and
$$
2\tau_2(m)=(\delta_1(1)-\delta_4(1))m-m(\mu_1(1)-\mu_4(1)),
$$
which are the required statements (3) and (5).

Let us choose $G=\left[\smallmatrix 1 & 0\\
n & 0
\endsmallmatrix \right]$ (resp. $G=\left[\smallmatrix 0 & 0\\
n & 1
\endsmallmatrix \right]$) in (3.2) and repeat the previous computational
process $(3.15)-(3.18)$. Then $\tau_3(n)=0$ will follow. We also
get
$$
\nu_3(n)=n\delta_1(1)+n\delta_3(n)-\mu_1(1)n-\mu_3(n)n  \eqno(3.19)
$$
and
$$
\nu_3(n)=\mu_3(n)n+\mu_4(1)n-n\delta_3(n)-n\delta_4(1). \eqno(3.20)
$$
Combining (3.19) with (3.20) gives
$$
n(\delta_1(1)+\delta_4(1)+2\delta_3(n))=(\mu_1(1)+\mu_4(1)+2\mu_3(n))n
$$
and
$$
2\nu_3(n)=n(\delta_1(1)-\delta_4(1))-(\mu_1(1)-\mu_4(1))n,
$$
which are the required statements (4) and (6).

Taking $G=\left[\smallmatrix a & m\\
0 & 0
\endsmallmatrix \right]$ into (3.2) yields
$$
[\delta_1(a), a]_k+[\delta_2(m), a]_k=0. \eqno(3.21)
$$
Since $\delta_1$ is a $k$-commuting mapping of $A$, $\delta_2(m)\in
{\mathcal Z}(A)_k$. Choosing $G=\left[\smallmatrix a & 0\\
n & 0
\endsmallmatrix \right]$ in (3.2) and using the same computational
methods, we can obtain $\delta_3(n)\in {\mathcal Z}(A)_k$. Likewise,
$\mu_2(m)\in {\mathcal Z}(B)_k$ and $\mu_3(n)\in {\mathcal Z}(B)_k$
will follow if we take $G=\left[\smallmatrix 0 & m\\
0 & b
\endsmallmatrix \right]$ and $G=\left[\smallmatrix 0 & 0\\
n & b
\endsmallmatrix \right]$ into (3.2), respectively. This completes the
proof of this proposition.
\end{proof}

The current authors in \cite{LiWei} described the general form of
arbitrary derivation on the generalized matrix algebra $\mathcal{G}$
and showed that every semi-centralizing derivation on $\mathcal{G}$
is zero. We next extend this result to the case of $k$-commuting
derivations.

\begin{proposition}{\rm \cite[Proposition 4.2]{LiWei}}\label{xxsec3.3}
An $\mathcal{R}$-linear mapping $\Theta_{\rm d}$ is a derivation of
$\mathcal{G}$ if and only if $\Theta_{\rm d}$ has the form
$$
\begin{aligned}
& \Theta_{\rm d}\left(\left[
\begin{array}
[c]{cc}%
a & m\\
n & b\\
\end{array}
\right]\right) \\
=& \left[
\begin{array}
[c]{cc}%
\delta_1(a)-mn_0-m_0n & am_0-m_0b+\tau_2(m)\\
n_0a-bn_0+\nu_3(n) & n_0m+nm_0+\mu_4(b)\\
\end{array}
\right] ,\\
& \forall \left[
\begin{array}
[c]{cc}%
a & m\\
n & b\\
\end{array}
\right]\in \mathcal{G},
\end{aligned} \eqno(\spadesuit)
$$
where $m_0\in M, n_0\in N$ and
$$
\begin{aligned} \delta_1:& A \longrightarrow A, &
 \tau_2: & M\longrightarrow M, & \nu_3: & N\longrightarrow N , &
\mu_4: & B\longrightarrow B
\end{aligned}
$$
are all $\mathcal{R}$-linear mappings satisfying the following
conditions:
\begin{enumerate}
\item[{\rm(1)}] $\delta_1$ is a derivation of $A$ with
$\delta_1(mn)=\tau_2(m)n+m\nu_3(n);$

\item[{\rm(2)}] $\mu_4$ is a derivation of $B$ with
$\mu_4(nm)=n\tau_2(m)+\nu_3(n)m;$

\item[{\rm(3)}] $\tau_2(am)=a\tau_{2}(m)+\delta_1(a)m$ and
$\tau_2(mb)=\tau_2(m)b+m\mu_4(b);$

\item[{\rm(4)}] $\nu_3(na)=\nu_3(n)a+n\delta_1(a)$ and
$\nu_3(bn)=b\nu_3(n)+\mu_4(b)n$
\end{enumerate}
\end{proposition}

\begin{proposition}\label{xxsec3.4}
Let $\mathcal{G}$ be a $2$-torsion free generalized matrix algebra.
Then every $k$-commuting derivation on $\mathcal {G}$ is zero.
\end{proposition}

\begin{proof}
Let $\Theta_{\rm d}$ be a $k$-commuting derivation on $\mathcal
{G}$. By Proposition \ref{xxsec3.3} we know that $\Theta_{\rm d}$ is
of the form
$$
\Theta_{\rm d}\left(\left[
\begin{array}
[c]{cc}%
a & m\\
n & b\\
\end{array}
\right]\right) = \left[
\begin{array}
[c]{cc}%
\delta_1(a)-mn_0-m_0n & am_0-m_0b+\tau_2(m)\\
n_0a-bn_0+\nu_3(n) & n_0m+nm_0+\mu_4(b)\\
\end{array}
\right] \eqno(3.22)
$$
for all $\left[
\begin{array}
[c]{cc}%
a & m\\
n & b\\
\end{array}
\right]\in \mathcal{G}$, where $m_0=\tau_1(1), n_0=\nu_1(1)$. Since
$\Theta_{\rm d}$ is $k$-commuting on $\mathcal{G}$,
$\tau_1(1)=\nu_1(1)=0$ by the relation (3.4). Therefore (3.22)
becomes
$$
\Theta_{\rm d}\left(\left[
\begin{array}
[c]{cc}%
a & m\\
n & b\\
\end{array}
\right]\right) = \left[
\begin{array}
[c]{cc}%
\delta_1(a) & \tau_2(m)\\
\nu_3(n) & \mu_4(b)\\
\end{array}
\right] \eqno(3.23)
$$
It should be remarked that $\delta_1$ and $\mu_4$ are derivations of
$A$ and $B$, respectively. Thus $\delta_1(1)=\mu_4(1)=0$. In view of
the conditions (5) and (6) in Proposition \ref{xxsec3.2} we know
that $\tau_2(m)=0$ for all $m\in M$ and $\nu_3(n)=0$ for all $n\in
N$. Furthermore, applying the condition (3) in Proposition
\ref{xxsec3.3} gives that $\delta_1(a)m=0$ for all $a\in A$ and
$m\in M$. Since $M$ is a faithful left $A$-module, $\delta_1(a)=0$
for all $a\in A$. Similarly, by the condition (4) in Proposition
\ref{xxsec3.3} we can obtain $\mu_4(b)=0$ for all $b\in B$. Hence,
$\Theta_{\rm d}$ has the form
$$
\Theta_{\rm d}\left(\left[
\begin{array}
[c]{cc}%
a & m\\
n & b\\
\end{array}
\right]\right) = \left[
\begin{array}
[c]{cc}%
0 & 0\\
0 & 0\\
\end{array}
\right]
$$
for all $\left[
\begin{array}
[c]{cc}%
a & m\\
n & b\\
\end{array}
\right]\in \mathcal{G}$, which is the desired result.
\end{proof}

Now we are in a position to state the main theorem of this article.
This result will provide a sufficient condition which enables
arbitrary $k$-commuting mapping on the generalized matrix algebra
$\mathcal{G}$ to be proper.

\begin{theorem}\label{xxsec3.5}
Let $\mathcal{G}$ be a $2$-torsion free generalized matrix algebra
and $\Theta$ be a $k$-commuting mapping of $\mathcal{G}$. If the
following three conditions are satisfied:
\begin{enumerate}
\item[{\rm(1)}] ${\mathcal Z}(A)_k=\pi_A({\mathcal Z}(\mathcal {G}));$
\item[{\rm(2)}] ${\mathcal Z}(B)_k=\pi_B({\mathcal Z}(\mathcal {G}));$
\item[{\rm(3)}] There exist $m_0\in M$ and $n_0\in N$ such that
$$
{\mathcal Z}(\mathcal{G})=\left\{\left[
\begin{array}[c]{cc}%
a & 0\\
0 & b
\end{array}
\right] \vline a\in {\mathcal Z}(A),b\in {\mathcal Z}(B),am_0=m_0b,
n_0a=bn_0 \right\},
$$
\end{enumerate}
then $\Theta$ is proper. That is, $\Theta$ has the form
$$
\Theta\left(\left[
\begin{array}[c]{cc}%
a & m\\
n & b
\end{array}
\right]\right)=\lambda \left[
\begin{array}[c]{cc}%
a & m\\
n & b
\end{array}
\right]+\zeta\left(\left[
\begin{array}[c]{cc}%
a & m\\
n & b
\end{array}
\right]\right), \forall \left[
\begin{array}[c]{cc}%
a & m\\
n & b
\end{array}
\right]\in \mathcal{G},
$$
where $\lambda\in {\mathcal Z}(\mathcal{G})$ and $\zeta:
\mathcal{G}\longrightarrow {\mathcal Z}(\mathcal{G})$ is an
$\mathcal{R}$-linear mapping.
\end{theorem}

\begin{proof} By
Proposition \ref{xxsec3.2} we know that $\Theta$ has the form
$(\bigstar)$. We will complete the proof of this theorem via the
following six steps.

\smallskip

{\bf Step 1.} $\delta_2(m)m=m\mu_2(m)$ and $n\delta_3(n)=\mu_3(n)n$
for all $m\in M$ and $n\in N$.

\bigskip

We first claim that
$$ \left[
\begin{array}[c]{cc}%
\delta_1(1)+\delta_4(1) & 0\\
0 & \mu_1(1)+\mu_4(1)
\end{array}
\right] \in Z(\mathcal{G}).\eqno(3.24)
$$
Actually, we have obtained $\delta_1(1) \in {\mathcal Z}(A)_k$ in
the proof of Proposition \ref{xxsec3.2}. Applying the conditions
(1) and (3) yields that $\delta_1(1) \in {\mathcal Z}(A)$.
Likewise, we also have $\delta_4(1)\in {\mathcal Z}(A)$ and hence
$\delta_1(1)+\delta_4(1) \in {\mathcal Z}(A)$. Similarly, we get
$\mu_1(1)+\mu_4(1)\in {\mathcal Z}(B)$. Moreover, it follows from
the statement (3) of Proposition \ref{xxsec3.2} that for arbitrary
element $m\in M$,
$$
\begin{aligned}
&(\delta_1(1)+\delta_4(1)+2\delta_2(m_0+m))(m_0+m)\\
&=(\delta_1(1)+\delta_4(1)+2\delta_2(m_0))m_0+2\delta_2(m)m_0\\
&\quad+2\delta_2(m_0)m+(\delta_1(1)+\delta_4(1)+2\delta_2(m))m\\
&=m_0(\mu_1(1)+\mu_4(1)+2\mu_2(m_0))+2\delta_2(m)m_0\\
&\quad+2\delta_2(m_0)m+m(\mu_1(1)+\mu_4(1)+2\mu_2(m)).
\end{aligned}\eqno(3.25)
$$
On the other hand, the statement (3) of Proposition \ref{xxsec3.2}
gives
$$
\begin{aligned}
&(\delta_1(1)+\delta_4(1)+2\delta_2(m_0+m))(m_0+m)\\
&=(m_0+m)(\mu_1(1)+\mu_4(1)+2\mu_2(m_0+m))\\
&=m_0(\mu_1(1)+\mu_4(1)+2\mu_2(m_0))+2m_0\mu_2(m)\\
&\quad+2m\mu_2(m_0)+m(\mu_1(1)+\mu_4(1)+2\mu_2(m)).
\end{aligned}\eqno(3.26)
$$
The above two equalities (3.25) and (3.26) imply that
$$
2\delta_2(m)m_0+2\delta_2(m_0)m=2m_0\mu_2(m)+2m\mu_2(m_0).\eqno(3.27)
$$
Taking $m=m_0$ into (3.27) leads to
$4\delta_2(m_0)m_0=4m_0\mu_2(m_0)$. Since $\mathcal{G}$ is
$2$-torsion free,
$$
\delta_2(m_0)m_0=m_0\mu_2(m_0).
$$
Thus the statement (3) of Proposition \ref{xxsec3.2} becomes
$$
(\delta_1(1)+\delta_4(1))m_0=m_0(\mu_1(1)+\mu_4(1)).\eqno(3.28)
$$
Likewise, we by the statement (4) of Proposition \ref{xxsec3.2}
arrive at
$$
2n_0\delta_3(n)+2n\delta_3(n_0)=2\mu_3(n)n_0+2\mu_3(n_0)n.\eqno(3.29)
$$
Thus $n_0\delta_3(n_0)=\mu_3(n_0)n_0$ will follow if we choose
$n=n_0$ in (3.29) and consider the $2$-torsion free property of
$\mathcal{G}$. Now the statement (4) of Proposition \ref{xxsec3.2}
becomes
$$
n_0(\delta_1(1)+\delta_4(1))=(\mu_1(1)+\mu_4(1))n_0.\eqno(3.30)
$$
Combining (3.28), (3.30) with condition (3) completes the proof of
(3.24). In view of the conditions (2), (3) and (3.24) we obtain
$$
\delta_2(m)m=m\mu_2(m)\hspace{10pt} {\rm and}\hspace{10pt}
n\delta_3(n)=\mu_3(n)n.\eqno(3.31)
$$
By the relations (3.17), (3.18) and (3.31) we get
$$
\tau_2(m)=\delta_1(1)m-m\mu_1(1)=m\mu_4(1)-\delta_4(1)m.\eqno(3.32)
$$
In view of the relations (3.20), (3.21) and (3.31) we have
$$
\nu_3(n)=n\delta_1(1)-\mu_1(1)n=\mu_4(1)n-n\delta_4(1).\eqno(3.33)
$$

\smallskip

{\bf Step 2.} $\delta_3(n)m=m\mu_3(n)$, $\mu_2(m)n=n\delta_2(m)$ for
all $m\in M$ and $n\in N$.

\bigskip

We assert that
$$
(\delta_3(n)m-m\mu_3(n))n=0\eqno(3.34)
$$
for all $m\in M, n\in N$. Proposition \ref{xxsec3.2} shows that
$\delta_3(n)\in {\mathcal Z}(A)_k$ for all $n\in N$. The conditions
(1) and (3) force that $\delta_3(n)\in {\mathcal Z}(A)$ for all
$n\in N$. Thus
\begin{align*}
&(\delta_3(n)m-m\mu_3(n))n=\delta_3(n)mn-m\mu_3(n)n\\
&=mn\delta_3(n)-m\mu_3(n)n=m(n\delta_3(n)-\mu_3(n)n)
\end{align*}
and the assertion follows from (3.31). Using the the same
computational method we conclude
$$
n(\delta_3(n)m-m\mu_3(n))=0,\eqno(3.35)
$$
$$
m(\mu_2(m)n-n\delta_2(m))=0\eqno(3.36)
$$
and
$$
(\mu_2(m)n-n\delta_2(m))m=0.\eqno(3.37)
$$

Let us choose $G=\left[\smallmatrix 1 & m\\
n & 0
\endsmallmatrix \right]$. It follows from (3.32), (3.33) and Proposition
\ref{xxsec3.2} that
$$
\Theta(G) =\left[
\begin{array}
[c]{cc}%
\delta_1(1)+\delta_2(m)+\delta_3(n) & \delta_1(1)m-m\mu_1(1) \\
n\delta_1(1)-\mu_1(1)n & \mu_1(1)+\mu_2(m)+\mu_3(n) \\
\end{array}
\right].
$$
Applying the relation (3.31) yields
$$
[\Theta(G), G] =\left[
\begin{array}
[c]{cc}%
\delta_1(1)mn-mn\delta_1(1) & \delta_3(n)m-m\mu_3(n) \\
\mu_2(m)n-n\delta_2(m) & nm\mu_1(1)-\mu_1(1)nm \\
\end{array}
\right].\eqno(3.38)
$$
It should be remarked that $\delta_1(1)\in {\mathcal Z}(A),
\mu_1(1)\in {\mathcal Z}(B), mn\in A$ and $nm\in B$. Thus (3.38)
becomes
$$
[\Theta(G), G] =\left[
\begin{array}
[c]{cc}%
0 & \delta_3(n)m-m\mu_3(n) \\
\mu_2(m)n-n\delta_2(m) & 0 \\
\end{array}
\right].\eqno(3.39)
$$
Consequently, for each $k\geq 2$, the relations (3.34)-(3.37)
jointly give
$$
\left[
\begin{array}
[c]{cc}%
0 & 0 \\
0 & 0 \\
\end{array}
\right]=[\Theta(G), G]_k\\
=\left[
\begin{array}
[c]{cc}%
0 & (-1)^{k+1}(\delta_3(n)m-m\mu_3(n)) \\
\mu_2(m)n-n\delta_2(m) & 0 \\
\end{array}
\right].
$$
Therefore
$$
\delta_3(n)m=m\mu_3(n) \quad \mu_2(m)n=n\delta_2(m).\eqno(3.40)
$$

\smallskip

{\bf Step 3.} $\left[
\begin{array}
[c]{cc}%
\delta_2(m) & 0 \\
0 & \mu_2(m) \\
\end{array}
\right]\in \mathcal{Z(G)}$ for all $m\in M$ and $\left[
\begin{array}
[c]{cc}%
\delta_3(n) & 0 \\
0 & \mu_3(n) \\
\end{array}
\right]\in \mathcal{Z(G)}$ for all $n\in N$.

\bigskip

It is not difficult to verify that
 $$\left[
\begin{array}
[c]{cc}%
\delta_2(m_0) & 0 \\
0 & \mu_2(m_0) \\
\end{array}
\right]\in \mathcal{Z(G)}.\eqno(3.41)
$$
Indeed, by Proposition \ref{xxsec3.2} and the conditions (1), (2)
and (3) we have
$$
\delta_2(m)\in {\mathcal Z}(A) \hspace{4pt} {\rm and} \hspace{4pt}
\mu_2(m)\in {\mathcal Z}(B)\eqno(3.42)
$$
for all $m\in M$. Furthermore, (3.31) and (3.40) imply that
$$
\delta_2(m_0)m_0=m_0\mu_2(m_0)\hspace{8pt} {\rm and} \hspace{8pt}
n_0\delta_2(m_0)=\mu_2(m_0)n_0.\eqno(3.43)
$$
Then (3.41) follows from (3.42), (3.43) and the condition (3).

By the definition of $\mathcal{Z(G)}$ and the relation (3.41) we
know that $\delta_2(m_0)m=m\mu_2(m_0)$. This forces (3.27) to be
$$
\delta_2(m)m_0=m_0\mu_2(m).\eqno(3.44)
$$
Combining the relations (3.40), (3.44) and the condition (3) leads
to
$$
\left[
\begin{array}
[c]{cc}%
\delta_2(m) & 0 \\
0 & \mu_2(m) \\
\end{array}
\right]\in \mathcal{Z(G)} \eqno(3.45)
$$
for all $m\in M$. Similarly, it can be proved that
$$
\left[
\begin{array}
[c]{cc}%
\delta_3(n) & 0 \\
0 & \mu_3(n) \\
\end{array}
\right]\in \mathcal{Z(G)} \eqno(3.46)
$$
for all $n\in N$.

\smallskip

{\bf Step 4.} For all $a\in A$, $m\in M$ and $n\in N$, we have
\begin{enumerate}
\item[{\rm(a)}] $\delta_1(a)m-m\mu_1(a)=a(\delta_1(1)m-m\mu_1(1))=a(m\mu_4(1)-\delta_4(1)m),$

\item[{\rm(b)}] $n\delta_1(a)-\mu_1(a)n=(n\delta_1(1)-\mu_1(1)n)a=(\mu_4(1)n-n\delta_4(1))a.$
\end{enumerate}

\bigskip

Let us choose $G=\left[\smallmatrix a & m\\
0 & 0
\endsmallmatrix \right]$. Then (3.31) and $\delta_2(m)\in {\mathcal Z}(A)$ imply
that
$$
[\Theta(G), G]=\left[
\begin{array}
[c]{cc}%
[\delta_1(a), a] & \alpha_1 \\
0 & 0 \\
\end{array}
\right],
$$
where $\alpha_1=\delta_1(a)m-m\mu_1(a)-a(\delta_1(1)m-m\mu_1(1))$.
It should be remarked that there exists a unique algebraic
isomorphism $\varphi: \pi_A(\mathcal{Z(G)})\longrightarrow
\pi_B(\mathcal{Z(G)})$ such that $am=m\varphi(a)$ and
$na=\varphi(a)n$ for all $a\in \pi_A(\mathcal{Z(G)}), m\in M$ and
$n\in N$. Therefore
$$
\alpha_1=\delta_1(a)m-\varphi^{-1}(\mu_1(a))m-a\delta_1(1)m+a\varphi^{-1}(\mu_1(1))m.
\eqno(3.47)
$$
On the other hand,  we by an inductive computation have
$$
[\Theta(G), G]_i=\left[
\begin{array}
[c]{cc}%
[\delta_1(a), a]_i & \alpha_i \\
0 & 0 \\
\end{array}
\right]\hspace{4pt} {\rm for} \hspace{4pt} {\rm each}
\hspace{4pt}i>1,
$$
where $\alpha_i=[\delta_1(a), a]_{i-1}m-a\alpha_{i-1}$. Then
$$
\begin{aligned}0 &=\alpha_k\\
&=[\delta_1(a), a]_{k-1}m-a[\delta_1(a),
a]_{k-2}m+\cdots+(-a)^{k-2}[\delta_1(a),
a]m+(-a)^{k-1}\alpha_1.
\end{aligned}
\eqno(3.48)$$ Combining (3.47) with (3.48) yields that
$$
\begin{aligned}
0&=([\delta_1(a), a]_{k-1}+\cdots+(-a)^{k-2}[\delta_1(a),
a]\\
&\hspace{10pt}+(-a)^{k-1}(\delta_1(a)-\varphi^{-1}(\mu_1(a))-a\delta_1(1)+a\varphi^{-1}(\mu_1(1))))m
\end{aligned}
$$
for all $m\in M$. The fact that $M$ is a faithful left $A$-module
leads to
\begin{align*}
&(-a)^{k-1}(\mu_1(a)+a\delta_1(1)-a\mu_1(1))\\
&=[\delta_1(a), a]_{k-1}+\cdots+(-a)^{k-2}[\delta_1(a),
a]+(-a)^{k-1}\delta_1(a).
\end{align*}
Then Proposition \ref{xxsec3.2} and the condition (1) show that
$(-a)^{k-1}(\mu_1(a)+a\delta_1(1)-a\mu_1(1))\in {\mathcal Z}(A)$.
That is,
$$
[\delta_1(a), a]_{k-1}+\cdots+(-a)^{k-2}[\delta_1(a),
a]+(-a)^{k-1}\delta_1(a)\in {\mathcal Z}(A).\eqno(3.49)
$$
We assert that $[\delta_1(a), a]=0$ for all $a\in A$. Indeed, it
follows from (3.49) that
$$
[[\delta_1(a), a]_{k-1}+\cdots+(-a)^{k-2}[\delta_1(a),
a]+(-a)^{k-1}\delta_1(a), a]_{k-1}=0 \eqno(3.50)
$$
for all $a\in A$. Since $\delta_1$ is a $k$-commuting on $A$,
$$
a^{k-1}[\delta_1(a), a]_{k-1}=0
$$
for all $a\in A$. Furthermore, the fact that $\delta_1(1)\in
{\mathcal Z}(A)$ gives $[\delta_1(a+1), a+1]_{k-1}=[\delta_1(a),
a]_{k-1}$ for all $a\in A$. By Lemma \ref{xxsec3.1} we get
$[\delta_1(a), a]_{k-1}=0$ for all $a\in A$. Repeating the same
process we arrive at
$$
0=[\delta_1(a), a]_{k-1}=[\delta_1(a), a]_{k-2}=\cdots=[\delta_1(a),
a].\eqno(3.51)
$$
Taking the relation $[\delta_1(a), a]=0$ into (3.48) we obtain
$a^{k-1}\alpha_1=0$. Let us fix arbitrary element $m\in M$ and
define $\alpha(a)=\delta_1(a)m-m\mu_1(a)-a(\delta_1(1)m-m\mu_1(1))$.
Then
$$
\begin{aligned}
\alpha(a+1)&=\delta_1(a+1)m-m\mu_1(a+1)-(a+1)(\delta_1(1)m-m\mu_1(1))\\
&=\delta_1(a)m-m\mu_1(a)+\delta_1(1)m-m\mu_1(1)\\
&\hspace{10pt}-a(\delta_1(1)m-m\mu_1(1))-\delta_1(1)m+m\mu_1(1)\\
&=\delta_1(a)m-m\mu_1(a)-a(\delta_1(1)m-m\mu_1(1))=\alpha(a)
\end{aligned}
$$
Applying Lemma \ref{xxsec3.1} again we have $\alpha(a)=0$. Thus
$$
\delta_1(a)m-m\mu_1(a)=a(\delta_1(1)m-m\mu_1(1)).\eqno(3.52)
$$
Combining (3.32) with (3.52) yields
$$
\delta_1(a)m-m\mu_1(a)=a(m\mu_4(1)-\delta_4(1)m).\eqno(3.53)
$$
This completes the proof of (a).
Likewise, if we take $G=\left[\smallmatrix a & 0\\
n & 0
\endsmallmatrix \right]$, then (3.31) and $\delta_3(n)\in {\mathcal Z}(A)$ imply
that
$$
[\Theta(G), G]=\left[
\begin{array}
[c]{cc}%
[\delta_1(a), a] & 0 \\
\beta_1 & 0 \\
\end{array}
\right],
$$
where $\beta_1=(n\delta_1(1)-\mu_1(1)n)a+\mu_1(a)n-n\delta_1(a)$. We
by an inductive computation conclude
$$
[\Theta(G), G]_i=\left[
\begin{array}
[c]{cc}%
[\delta_1(a), a]_i & 0 \\
\beta_i & 0 \\
\end{array}
\right]{\rm for} \hspace{4pt} {\rm each} \hspace{4pt}i>1,
$$
where $\beta_i=\beta_{i-1}a-n[\delta_1(a), a]_{i-1}$. Then
$$
0=\beta_k\\
=\beta_1a^{k-1}-n[\delta_1(a), a]a^{k-2}-\cdots-n[\delta_1(a),
a]_{k-1}. \eqno(3.54)
$$
Taking (3.51) into (3.54) we obtain $\beta_1a^{k-1}=0$. Let us fix
arbitrary element $n\in N$ and define
$\beta(a)=(n\delta_1(1)-\mu_1(1)n)a+\mu_1(a)n-n\delta_1(a)$. A
direct computations gives $\beta(a)=\beta(a+1)$. In view of Lemma
\ref{xxsec3.1}, we have $\beta=0$. That is,
$$
(n\delta_1(1)-\mu_1(1)n)a=n\delta_1(a)-\mu_1(a)n. \eqno(3.55)
$$
It follows from the relations (3.33) and (3.55) that
$$
(\mu_4(1)n-n\delta_4(1))a=n\delta_1(a)-\mu_1(a)n, \eqno(3.56)
$$
which is the desired result (b).

\smallskip

{\bf Step 5.} For all $b\in B$, $m\in M$ and $n\in N$, we have

\begin{enumerate}
\item[{\rm(c)}]$\delta_4(b)m-m\mu_4(b)=(m\mu_1(1)-\delta_1(1)m)b=(\delta_4(1)m-m\mu_4(1))b,$

\item[{\rm(d)}]$n\delta_4(b)-\mu_4(b)n=b(\mu_1(1)n-n\delta_1(1))=b(n\delta_4(1)-\mu_4(1)n).$
\end{enumerate}

\bigskip

Let us choose $G=\left[\smallmatrix 0 & m\\
0 & b
\endsmallmatrix \right]$. Then (3.31) and $\mu_2(m)\in {\mathcal Z}(A)$ imply
that
$$
[\Theta(G), G]=\left[
\begin{array}
[c]{cc}%
0 & \gamma_1 \\
0 & [\mu_4(b), b] \\
\end{array}
\right],
$$
where $\gamma_1=\delta_4(b)m-m\mu_4(b)-(m\mu_1(1)-\delta_1(1)m)b$.
Note that there exists a unique algebraic isomorphism $\varphi:
\pi_A(\mathcal{Z(G)})\longrightarrow \pi_B(\mathcal{Z(G)})$ such
that $am=m\varphi(a)$ and $na=\varphi(a)n$ for all $a\in
\pi_A(\mathcal{Z(G)}), m\in M$ and $n\in N$. Therefore
$$
\gamma_1=m\varphi(\delta_4(b))-m\mu_4(b)-m\mu_1(1)b-m\varphi(\delta_1(1))b.\eqno(3.57)
$$
On the other hand,  we by an inductive computation get
$$
[\Theta(G), G]_i=\left[
\begin{array}
[c]{cc}%
0 & \gamma_i \\
0 & [\mu_4(b), b]_i \\
\end{array}
\right] {\rm for} \hspace{4pt} {\rm each} \hspace{4pt}i>1,
$$
where $\gamma_i=\gamma_{i-1}b-m[\mu_4(b), b]_{i-1}$. Then
$$
0=\gamma_k=\gamma_1b^{k-1}-m[\mu_4(b), b]b^{k-2}-\cdots-m[\mu_4(b),
b]_{k-1}.\eqno(3.58)
$$
Combining (3.57) with (3.58) yields that
$$
\begin{aligned}
0&=m((\varphi(\delta_4(b))-\mu_4(b)-\mu_1(1)b-\varphi(\delta_1(1)b))b^{k-1}
\\
&\hspace{10pt}-[\mu_4(b), b]b^{k-2}-\cdots-[\mu_4(b), b]_{k-1})
\end{aligned}
$$
for all $m\in M$. Since $M$ is a faithful right $B$-module, we have
\begin{align*}
&(\varphi(\delta_4(b))-\mu_1(1)b-\varphi(\delta_1(1)b))b^{k-1}
\\&=\mu_4(b)b^{k-1}+[\mu_4(b), b]b^{k-2}+\cdots+[\mu_4(b), b]_{k-1}.
\end{align*}
Then Proposition \ref{xxsec3.2} and the conditions (1) and (2)
jointly give
$(\varphi(\delta_4(b))-\mu_1(1)b-\varphi(\delta_1(1)b)b^{k-1}\in
{\mathcal Z}(B)$. That is,
$$
\mu_4(b)b^{k-1}+[\mu_4(b), b]b^{k-2}+\cdots+[\mu_4(b), b]_{k-1}\in
{\mathcal Z}(B).\eqno(3.59)
$$
We conclude that $[\mu_4(b), b]=0$ for all $b\in B$. Indeed, it
follows from the relation (3.59) that
$$
[\mu_4(b)b^{k-1}+[\mu_4(b), b]b^{k-2}+\cdots+[\mu_4(b), b]_{k-1},
b]_{k-1}=0\eqno(3.60)
$$
for all $b\in B$. Since $\mu_4$ is $k$-commuting on $B$,
$$
[\mu_4(b), b]_{k-1}b^{k-1}=0
$$
for all $b\in B$. Furthermore, $\mu_4(1)\in {\mathcal Z}(B)$ leads
to $[\mu_4(b+1), b+1]_{k-1}=[\mu_4(b), b]_{k-1}$ for all $b\in B$.
Applying Lemma \ref{xxsec3.1} we obtain $[\mu_4(b), b]_{k-1}=0.$
Repeating the same process as above we arrive at
$$
0=[\mu_4(b), b]_{k-1}=[\mu_4(b), b]_{k-2}=\cdots=[\mu_4(b),
b].\eqno(3.61)
$$
Taking $[\mu_4(b), b]=0$ into (3.58) we get $\gamma_1b^{k-1}=0$. Let
us fix arbitrary element $m\in M$ and define
$\gamma(b)=\delta_4(b)m-m\mu_4(b)-(m\mu_1(1)-\delta_1(1)m)b$. A
direct computations gives $\gamma(b+1)=\gamma(b)$. And hence
$\gamma(b)=0$ by Lemma \ref{xxsec3.1}. That is,
$$
\delta_4(b)m-m\mu_4(b)=(m\mu_1(1)-\delta_1(1)m)b.\eqno(3.62)
$$
It follows from the relations (3.32) and (3.62) that
$$
\delta_4(b)m-m\mu_4(b)=(\delta_4(1)m-m\mu_4(1))b.\eqno(3.63)
$$
This completes the proof of (c).

In order to prove (d), let us take $G=\left[\smallmatrix 0 & 0\\
n & b
\endsmallmatrix \right]$. Then (3.31) and $\mu_3(n)\in {\mathcal Z}(B)$ imply
that
$$
[\Theta(G), G]=\left[
\begin{array}
[c]{cc}%
0 & 0 \\
\eta_1 & [\mu_4(b), b] \\
\end{array}
\right],
$$
where $\eta_1=\mu_4(b)n-n\delta_4(b)-b(\mu_1(1)n-n\delta_1(1))$. It
is easy to verify that for each $i>1$,
$$
[\Theta(G), G]_i=\left[
\begin{array}
[c]{cc}%
0 & 0 \\
\eta_i & [\mu_4(b), b]_i \\
\end{array}
\right],
$$
where $\eta_i=[\mu_4(b), b]_{i-1}n-b\eta_{i-1}$. Therefore
$$
\begin{aligned}
0&=\eta_k\\
&=[\mu_4(b), b]_{k-1}n-b[\mu_4(b),
b]_{k-2}n+\cdots+(-b)^{k-2}[\mu_4(b), b]n+(-b)^{k-1}\eta_1.
\end{aligned} \eqno(3.64)
$$
In view of the relations (3.61) and (3.64) we have
$b^{k-1}\eta_1=0$. Let us fix arbitrary element $n\in N$ and define
$\eta(b)=\mu_4(b)n-n\delta_4(b)-b(\mu_1(1)n-n\delta_1(1))$. Then a
straightforward computations gives $\eta(b)=\eta(b+1)$. By Lemma
\ref{xxsec3.1} again it follows that $\eta=0$. That is,
$$
\mu_4(b)n-n\delta_4(b)=b(\mu_1(1)n-n\delta_1(1)).\eqno(3.65)
$$
Combining (3.33) with (3.65) leads to
$$
n\delta_4(b)-\mu_4(b)n=b(n\delta_4(1)-\mu_4(1)n).\eqno(3.66)
$$
This completes the proof of this step.

\smallskip

{\bf Step 6.} $\Theta$ is proper.

\bigskip

Suppose that
$$
\Omega(X):=\Theta(X)-XC
$$
for all $X\in \mathcal{G}$, where $C=\left[\smallmatrix
\delta_1(1)-\varphi^{-1}(\mu_1(1)) & 0\\
 0 & \varphi(\delta_1(1))-\mu_1(1)
\endsmallmatrix
\right]\in \mathcal{Z(G)}$. We assert that
$\Omega(\mathcal{G})\subseteq \mathcal{Z(G)}$. Note that there
exists a unique algebraic isomorphism $\varphi:
\pi_A(\mathcal{Z(G)})\longrightarrow \pi_B(\mathcal{Z(G)})$ such
that $am=m\varphi(a)$ and $na=\varphi(a)n$ for all $a\in
\pi_A(\mathcal{Z(G)}), m\in M$ and $n\in N$. Thus we get
$$
\begin{aligned}
\Omega\left(\left[\begin{array}[c]{cc}%
a & m\\
n & b
\end{array}
\right]\right)=& \left[\begin{array}[c]{cc}%
\delta_1(a)-a\delta_1(1)+a\varphi^{-1}(\mu_1(1)) & 0\\
0 & \mu_1(a)
\end{array}
\right] \\
&+ \left[\begin{array}[c]{cc}%
\delta_4(b) & 0\\
0 & \mu_4(b)-b\varphi(\delta_1(1))+b\mu_1(1)
\end{array}
\right] \\
&+ \left[\begin{array}[c]{cc}%
\delta_2(m)+\delta_3(n) & 0\\
0 & \mu_2(m)+\mu_3(n)
\end{array}
\right]
\end{aligned} \eqno(3.67)
$$
for all $\left[\smallmatrix
a & m\\
n & b
\endsmallmatrix
\right]\in \mathcal{G}$. By the relation (3.52) we know that
$$
\begin{aligned}
&\hspace{12pt}(\delta_1(a)-a\delta_1(1)+a\varphi^{-1}(\mu_1(1)))m-m\mu_1(a)\\
&=a(\delta_1(1)m-m\mu_1(1))-a(\delta_1(1)-\varphi^{-1}(\mu_1(1)))m\\
&=0
\end{aligned}
$$
for all $a\in A, m\in M$. That is,
$$
(\delta_1(a)-a\delta_1(1)+a\varphi^{-1}(\mu_1(1)))m=m\mu_1(a)
\eqno(3.68)
$$
for all $a\in A, m\in M$. In view of the relation (3.55) and the
fact $\delta_1(1)\in {\mathcal Z}(A)$ we obtain
$$
\begin{aligned}
&\hspace{12pt}n(\delta_1(a)-a\delta_1(1)+a\varphi^{-1}(\mu_1(1)))-\mu_1(a)n\\
&=(n\delta_1(1)-\mu_1(1)n)a-na(\delta_1(1)-\varphi^{-1}(\mu_1(1)))\\
&=0
\end{aligned}
$$
for all $a\in A, n\in N$. Therefore
$$
n(\delta_1(a)-a\delta_1(1)+a\varphi^{-1}(\mu_1(1)))=\mu_1(a)n
\eqno(3.69)
$$
for all $a\in A, n\in N$. In view of (3.68), (3.69) and the
definition of $\mathcal{Z(G)}$ we conclude
$$
\left[\begin{array}[c]{cc}%
\delta_1(a)-a\delta_1(1)+a\varphi^{-1}(\mu_1(1)) & 0\\
0 & \mu_1(a)
\end{array}
\right]\in \mathcal{Z(G)} \eqno(3.70)
$$
for all $a\in A$. Likewise, the relation (3.62) and the fact
$\mu_1(1)\in {\mathcal Z}(B)$ jointly lead to
$$
\begin{aligned}
&\hspace{12pt}\delta_4(b)m-m(\mu_4(b)-b\varphi(\delta_1(1))+b\mu_1(1))\\
&=(m\mu_1(1)-\delta_1(1)m)b+mb(\varphi(\delta_1(1))-\mu_1(1))\\
&=0
\end{aligned}
$$
for all $b\in B, m\in M$. This implies that
$$
\delta_4(b)m=m(\mu_4(b)-b\varphi(\delta_1(1))+b\mu_1(1))\eqno(3.71)
$$
for all $b\in B, m\in M$. It follows from (3.65) that
$$
\begin{aligned}
&\hspace{12pt}(\mu_4(b)-b\varphi(\delta_1(1))+b\mu_1(1))n-n\delta_4(b)\\
&=b(n\delta_1(1)-\mu_1(1)n)-b(\varphi(\delta_1(1))-\mu_1(1))n\\
&=0
\end{aligned}
$$
for all $b\in B, n\in N$. This shows
$$
n\delta_4(b)=(\mu_4(b)-b\varphi(\delta_1(1))+b\mu_1(1))n \eqno(3.72)
$$
for all $b\in B, n\in N$. Taking into account (3.71), (3.72) and the
definition of $\mathcal{Z(G)}$ yields
$$
\left[\begin{array}[c]{cc}%
\delta_4(b) & 0\\
0 & \mu_4(b)-b\varphi(\delta_1(1))+b\mu_1(1)
\end{array}
\right]\in \mathcal{Z(G)} \eqno(3.73)
$$
for all $b\in B$. Then (3.45), (3.46), (3.70) and (3.73) prove that
$$
\Omega\left(\left[\begin{array}[c]{cc}%
a & m\\
n & b
\end{array}
\right]\right)\in \mathcal{Z(G)}
$$
for all $\left[\begin{array}[c]{cc}%
a & m\\
n & b
\end{array}
\right]\in \mathcal{G}$, which is the desired assertion.
\end{proof}

The following result will be used in the sequel.

\begin{corollary}\label{xxsec3.6}
If ${\mathcal Z}(A)_k=\mathcal{R}1={\mathcal Z}(B)_k$, then every
$k$-commuting mapping of $\mathcal{G}$ is proper.
\end{corollary}

Indeed, since ${\mathcal R}1\subseteq \pi_A(\mathcal{Z(G)})\subseteq
{\mathcal Z}(A)\subseteq {\mathcal Z}(A)_k={\mathcal R}1$, we have
${\mathcal R}1={\mathcal Z}(A)={\mathcal
Z}(A)_k=\pi_A(\mathcal{Z(G)})$. Likewise, we get ${\mathcal
R}1={\mathcal Z}(B)={\mathcal Z}(B)_k=\pi_B(\mathcal{Z(G)})$. It is
easy to check that the condition (3) of Theorem \ref{xxsec3.5} is
satisfied. This corollary follows from Theorem \ref{xxsec3.5}.
\vspace{2mm}

In particular, if the generalized matrix algebra $\mathcal{G}$
degenerates one general triangular algebra (that is, ${\mathcal
G}=\left[\begin{array}[c]{cc}%
A & M\\
O & B
\end{array}
\right]$), then our main result Theorem \ref{xxsec3.5} contains the
main theorem of \cite{DuWang1} as a special case.

\begin{corollary}{\rm \cite[Theorem 1.1]{DuWang1}}\label{xxsec3.7}
Let $\mathcal{T}$ be a $2$-torsion free triangular algebra and
$\Theta$ be a $k$-commuting mapping of $\mathcal{T}$. If the
following three conditions are satisfied:
\begin{enumerate}
\item[{\rm(1)}] ${\mathcal Z}(A)_k=\pi_A({\mathcal Z}(\mathcal {G}));$
\item[{\rm(2)}] ${\mathcal Z}(B)_k=\pi_B({\mathcal Z}(\mathcal {G}));$
\item[{\rm(3)}] There exists $m_0\in M$ such that
$$
{\mathcal Z}(\mathcal{G})=\left\{\left[
\begin{array}[c]{cc}%
a & 0\\
0 & b
\end{array}
\right] \vline a\in {\mathcal Z}(A),b\in {\mathcal Z}(B),am_0=m_0b
\right\},
$$
\end{enumerate}
then $\Theta$ is proper.
\end{corollary}

\bigskip

\section{Applications}\label{xxsec4}

In this section, we will present some applications of $k$-commuting
mappings to full matrix algebras, inflated algebras , upper and
lower triangular matrix algebras, nest algebras and block upper and
lower triangular matrix algebras.

\subsection{Full matrix algebras}
\label{xxsec4.1}

Let $\mathcal{R}$ be a commutative ring with identity, $A$ be a
$2$-torsion free unital algebra over $\mathcal{R}$ and $M_n(A)$ be
the algebra of $n\times n$ matrices with $n\geq 2$. Then the
\textit{full matrix algebra} $M_n(A)(n\geq 2)$ can be represented as
a generalized matrix algebra of the form
$$
M_n(A)=\left[
\begin{array}
[c]{cc}%
A & M_{1\times(n-1)}(A)\\
M_{(n-1)\times 1}(A) & M_{n-1}(A)\\
\end{array}
\right].
$$

\begin{corollary}\label{xxsec4.1}
Every $k$-commuting mapping on the full matrix algebra $M_n(A)$ or
$M_n(\mathcal{R})$ is proper.
\end{corollary}

One can directly check that $M_n(A)$ or $M_n(\mathcal{R})$ satisfies
all conditions (1)-(3) of Theorem \ref{xxsec3.5}. Therefore, every
commuting mapping on $M_n(A)$ or $M_n(\mathcal{R})$ is proper. We
would like to point out that this corollary can also be obtained by
applying the notion of ${\rm FI}$-degree of functional identities
and related results in \cite{BresarChebotarMartindale}.

\subsection{Inflated algebras}\label{xxsec4.2}

Let $A$ be a unital $\mathcal{R}$-algebra and $V$ be an
$\mathcal{R}$-linear space. Given an $\mathcal{R}$-bilinear form
$\gamma: V\otimes_{\mathcal{R}}V\rightarrow A$, we define an
associative algebra (not necessarily with identity)
$B=B(A,V,\gamma)$ as follows: As an $\mathcal{R}$-linear space,
$B$ equals to $V\otimes_{\mathcal{R}}V\otimes_{\mathcal{R}}A$. The
multiplication is defined as follows
$$
(a\otimes b\otimes x)\cdot(c\otimes d\otimes y):=a\otimes d\otimes
x\gamma(b,c)y
$$
for all $a,b,c,d\in V$ and any $x, y\in A$. This definition makes
$B$ become an associative $\mathcal{R}$-algebra and $B$ is called
an \textit{inflated algebra} of $A$ along $V$. The inflated
algebras are closely connected with the cellular algebras which
are extensively studied in representation theory. We refer the
reader to \cite{KonigXi} and the references therein for these
algebras.

Let us assume that $V$ is a non-zero linear space with a basis
$\{v_1,\cdots,v_n\}$. Then the bilinear form $\gamma$ can be
characterized by an $n\times n$ matrix $\Gamma$ over $A$, that is,
$\Gamma=(\gamma(v_i,v_j))$ for $1\leq i,j\leq n$. Now we could
define a new multiplication $``\circ"$ on the full matrix algebra
$M_n(A)$ by
$$
X\circ Y:=X\Gamma Y\quad\text{for all}\quad X,Y\in M_n(A).
$$
Under the usual matrix addition and the new multiplication
$``\circ"$, $M_n(A)$ becomes a new associative algebra which is a
generalized matrix algebra in the sense of Brown \cite{Brown}. We
denote this new algebra by $(M_n(A),\Gamma)$. It should be
remarked that our current generalized matrix algebras contain all
generalized matrix algebras defined by Brown \cite{Brown} as
special cases. By \cite[Lemma 4.1]{KonigXi}, the inflated algebra
$B(A,V,\gamma)$ is isomorphic to $(M_n(A),\Gamma)$ and hence is a
generalized matrix algebra in the sense of ours.

\begin{corollary}\label{xxsec4.2}
Let $A$ be a unital $\mathcal{R}$-algebra, $V$ be an
$\mathcal{R}$-linear space and $B(A,V,\gamma)$ be the inflated
algebra of $A$ along $V$. If $B(A,V,\gamma)$ has an identity
element, then each $k$-commuting mapping of $B(A,V,\gamma)$ is
proper.
\end{corollary}

\begin{proof}
If $B(A,V,\gamma)$ has an identity element, then the matrix
$\Gamma$ defined by the bilinear form $\gamma$ is invertible in
the full matrix algebra $M_n(A)$ by \cite[Proposition
4.2]{KonigXi}. We define
$$
\begin{aligned}
\sigma: M_n(A) & \longrightarrow (M_n(A),\Gamma)\\
X & \longmapsto X\Gamma^{-1}.
\end{aligned}
$$
Note that $\sigma(X)\circ \sigma(Y)=\sigma(X)\Gamma
\sigma(Y)=XY\Gamma^{-1}=\sigma(XY)$ for all $X,Y\in M_n(A)$ and
hence $\sigma$ is an algebraic isomorphism. Now the result follows
from Corollary \ref{xxsec4.1} and the fact $B(A,V,\gamma)\cong
(M_n(A),\Gamma)$.
\end{proof}

\subsection{Prime algebras}
\label{xxsec4.3}

Let $\mathcal{A}$ be a $2$-torsionfree prime algebra over a commutative ring $\mathcal{R}$. Assume that $\mathcal{A}$ has
unit $1$ and a nontrivial idempotent $e$. Then $\mathcal{A}=e\mathcal{A}e+e\mathcal{A}(1-e)+(1-e)\mathcal{A}e+(1-e)\mathcal{A}(1-e)$
Since $\mathcal{A}$ is prime, it is straightforward to verify that both $e\mathcal{A}e$ and $(1-e)\mathcal{A}(1-1)$ are also prime.
By \cite[Theorem 2]{Lanski} it follows that $\mathcal{Z}(e\mathcal{A}e)_k=\pi_{e\mathcal{A}e}\mathcal{Z}(\mathcal{A})$ and
$\mathcal{Z}((1-e)\mathcal{A}(1-e))_k=\pi_{(1-e)\mathcal{A}(1-e)}\mathcal{Z}(\mathcal{A})$.

\begin{corollary}\label{xxsec4.3}
Let $\mathcal{A}$ be a $2$-torsionfree prime algebra over a commutative ring $\mathcal{R}$. Suppose that $\mathcal{A}$ has
the identity $1$ and a nontrivial idempotent $e$. Then every $k$-commuting mapping of $\mathcal{A}$ is proper.
\end{corollary}

\subsection{Upper and lower matrix triangular algebras}
\label{xxsec4.4}

Let $\mathcal{R}$ be a $2$-torsion free commutative ring with
identity. We denote the set of all $p\times q$ matrices over
$\mathcal{R}$ by $M_{p\times q}(\mathcal{R})$. Let us denote the set
of all $n\times n$ upper triangular matrices over $\mathcal{R}$ and
the set of all $n\times n$ lower triangular matrices over
$\mathcal{R}$ by $\mathcal{T}_n(\mathcal{R})$ and
$\mathcal{T}_n^\prime(\mathcal{R})$, respectively. For $n\geq 2$ and
each $1\leq k \leq n-1$, the \textit{upper triangular matrix
algebra} $\mathcal{T}_n(\mathcal{R})$ and \textit{lower triangular
matrix algebra} $\mathcal{T}_n^\prime(\mathcal{R})$ can be written
as
$$
\mathcal{T}_n(\mathcal{R})=\left[
\begin{array}
[c]{cc}%
\mathcal{T}_k(\mathcal{R}) & M_{k\times (n-k)}(\mathcal{R})\\
 & \mathcal{T}_{n-k}(\mathcal{R})
\end{array}
\right] \hspace{5pt} {\rm and} \hspace{5pt}
\mathcal{T}_n^\prime(\mathcal{R})=\left[
\begin{array}
[c]{cc}%
\mathcal{T}_k^\prime(\mathcal{R}) &  \\
M_{(n-k)\times k}(\mathcal{R}) &
\mathcal{T}_{n-k}^\prime(\mathcal{R})
\end{array}
\right],
$$
respectively.

\begin{corollary}\label{xxsec4.4}
Every $k$-commuting mapping of the upper triangular matrix algebra
$\mathcal{T}_n(\mathcal{R})$ {\rm (} resp. the lower triangular
matrix algebra $\mathcal{T}_n^\prime(\mathcal{R})${\rm )} is proper.
\end{corollary}

We will give a unification proof for the cases of the upper and
lower triangular matrix algebras and nest algebras in below.

\subsection{Nest algebras}
\label{xxsec4.5}

Let $\mathbf{H}$ be a complex Hilbert space and
$\mathcal{B}(\mathbf{H})$ be the algebra of all bounded linear
operators on $\mathbf{H}$. Let $I$ be a index set. A \textit{nest}
is a set $\mathcal{N}$ of closed subspaces of $\mathbf{H}$
satisfying the following conditions:
\begin{enumerate}
\item[(1)] $0, \mathbf{H}\in \mathcal{N}$; \item[(2)] If $N_1,
N_2\in \mathcal{N}$, then either $N_1\subseteq N_2$ or
$N_2\subseteq N_1$; \item[(3)] If $\{N_i\}_{i\in I}\subseteq
\mathcal{N}$, then $\bigcap_{i\in I}N_i\in \mathcal{N}$;
\item[(4)] If $\{N_i\}_{i\in I}\subseteq \mathcal{N}$, then the
norm closure of the linear span of $\bigcup_{i\in I} N_i$ also
lies in $\mathcal{N}$.
\end{enumerate}
If $\mathcal{N}=\{0, \mathbf{H}\}$, then $\mathcal{N}$ is called a
trivial nest, otherwise it is called a non-trivial nest.

The \textit{nest algebra} associated with $\mathcal{N}$ is the set
$$
\mathcal{T}(\mathcal{N})=\{\hspace{3pt} T\in
\mathcal{B}(\mathbf{H})\hspace{3pt}| \hspace{3pt} T(N)\subseteq N
\hspace{3pt} {\rm for} \hspace{3pt} {\rm all} \hspace{3pt} N\in
\mathcal{N}\} .
$$
A nontrivial nest algebra is a triangular algebra. Indeed, if
$N\in \mathcal{N}\backslash \{0, {\mathbf H}\}$ and $E$ is the
orthogonal projection onto $N$, then
$\mathcal{N}_1=E(\mathcal{N})$ and
$\mathcal{N}_2=(1-E)(\mathcal{N})$ are nests of $N$ and
$N^{\perp}$, respectively. Moreover,
$\mathcal{T}(\mathcal{N}_1)=E\mathcal{T}(\mathcal{N})E,
\mathcal{T}(\mathcal{N}_2)=(1-E)\mathcal{T}(\mathcal{N})(1-E)$ are
nest algebras and
$$
\mathcal{T}(\mathcal{N})=\left[
\begin{array}
[c]{cc}%
\mathcal{T}(\mathcal{N}_1) & E\mathcal{T}(\mathcal{N})(1-E)\\
O & \mathcal{T}(\mathcal{N}_2)\\
\end{array}
\right].
$$
Note that any finite dimensional nest algebra is isomorphic to a
complex block upper triangular matrix algebra. We refer the reader
to \cite{Davidson} for the theory of nest algebras.

\begin{corollary}\label{xxsec4.5}
Every $k$-commuting mapping of the nest algebra
$\mathcal{T}(\mathcal{N})$ is proper.
\end{corollary}

We now give a unification proof for Corollary \ref{xxsec4.3} and
Corollary \ref{xxsec4.4} by an induction on $k$. For convenience,
let us set $\mathcal{W}=\mathcal{T}_n(\mathcal{R}),
\mathcal{T}_n^\prime(\mathcal{R})$ or $\mathcal{T}(\mathcal{N})$. The case
of $k=1$ is clearly trivial, since
$\mathcal{Z}(\mathcal{W})_1=\mathcal{R}1$. Let us choose an
arbitrary element $W\in \mathcal{Z(W)}_k$. Then $[[W, X],
X]_{k-1}=0$ for all $X\in \mathcal{W}$. If $\mathcal{W}$ is a
trivial nest algebra, then $\mathcal{W}={\mathcal B}(\mathbf{H})$ is
a centrally closed prime algebra. By \cite[Theorem 1]{Lanski} it
follows that $[W,X]=0$ for all $X\in \mathcal{W}$. This implies that
$W\in \mathcal{R}1$ and that $\mathcal{Z(W)}_k=\mathcal{R}1$. If
$\mathcal{W}$ is a nontrivial nest algebra or an upper triangular
matrix algebra. Then $\mathcal{W}$ can be written as the triangular
algebra
$$
\left[
\begin{array}
[c]{cc}%
A & M\\
O & B\\
\end{array}
\right].
$$
In view of the induction hypothesis we have
$\mathcal{Z}(A)_{k-1}=\mathcal{R}1=\mathcal{Z}(B)_{k-1}$. By
Corollary \ref{xxsec3.6} we know that there exist $\lambda\in
\mathcal{R}$ and $\zeta: \mathcal{W}\longrightarrow \mathcal{R}1$
such that
$$
[W, X]=\lambda X+\zeta(X)
$$
for all $X\in \mathcal{W}$. Therefore
$$
(W-\lambda I)X+X(-W)\in \mathcal{Z(W)}=\mathcal{R}1
$$
for all $X\in \mathcal{W}$. A straightforward computation leads to
$(W-\lambda I)=-(-W)\in \mathcal{Z(W)}$. This shows that
$\mathcal{Z(W)}_k=\mathcal{R}1$. Thus
$\mathcal{Z}(A)_k=\mathcal{R}1=\mathcal{Z}(B)_k$. It follows from
Corollary \ref{xxsec3.6} that every $k$-commuting mapping on the
upper (resp. lower) triangular matrix algebra
$\mathcal{T}_n(\mathcal{R})$ (resp.
$\mathcal{T}_n^\prime(\mathcal{R})$) is proper. Likewise, every
$k$-commuting mapping on the nest algebra $\tau(\mathcal{N})$ is
also proper.

\subsection{von Neumann Algebras}
\label{xxsec4.6}

Recall that a von Neumann algebra $M$ is a subalgebra of some $\mathcal{B}(\mathbf{H})$,
the algebra of all bounded linear operators acting on a complex Hilbert space $\mathbf{H}$,
which satisfies the double commutant property:  $\mathcal{M}^{\prime\prime}=\mathcal{M}$,
where $\mathcal{M}^\prime=\{ T\in \mathcal{B}(\mathbf{H}) |  TA=AT,  \forall A\in \mathcal{M}  \}$
and $\mathcal{M}^{\prime\prime}=\{\mathcal{M}^\prime\}^\prime$. For $A\in \mathcal{M}$,
the central carrier of $A$, denoted by $\overline{A}$, is the intersection of all central projections
$P$ such that $PA=0$. If $A$ is self-adjoint, then the core of $A$, denoted by $\underline{A}$,
is ${\rm sup}\{ S\in \mathcal{Z(M)}| S=S^\ast, S\leq A \}$. In case of $A=P$ is a projection, then
$\underline{P}$ is the largest central projection $\leq P$. A projection $P$ is core-free if $\underline{P}=0$.
It is easy to see that $\underline{P}=0$ if and only if $\overline{I-P}=I$. If $\mathcal{M}$ has no
central summands of type $I_1$, then each nonzero central projection in $\mathcal{M}$ is the
carrier of a core-free projection in $\mathcal{M}$. In particular, there exists a nonzero core-free
projection $P\in  \mathcal{M}$ with $\overline{P}=I$. For such $P$, note that $\overline{P}=\overline{I-P}=I$.
It follows from the definition of the central carrier that both ${\rm span}\{TP(x) | T\in \mathcal{M}, x\in \mathbf{H} \}$
and ${\rm span}\{ T(I -P)(x) | T\in \mathcal{M}, x\in \mathbf{H} \}$ are dense in $\mathbf{H}$. So $A\mathcal{M}P=0$
implies $A=0$ and $A\mathcal{M}(I-P)=0$ implies $A=0$. Thus, if $\mathcal{M}$ has no central summands
of type $I_1$, then $\mathcal{M}=P\mathcal{M}P+P\mathcal{M}Q+Q\mathcal{M}P+Q\mathcal{M}Q$
satisfies the corresponding conditions (1)-(3) in Theorem \ref{xxsec3.5}.
This is due to the fact $\mathcal{M}$ is a semiprime algebra and $P\mathcal{M}P$ and $Q\mathcal{M}Q$ are both semiprime.
Applying \cite[Theorem 2]{Lanski} yields that ${\mathcal Z}(P\mathcal{M}P)_k=\pi_{P\mathcal{M}P}({\mathcal Z}(\mathcal {M}))$
and ${\mathcal Z}(Q\mathcal{M}Q)_k=\pi_{Q\mathcal{M}Q}({\mathcal Z}(\mathcal {M}))$.
Therefore, Theorem \ref{xxsec3.5}  is true for additive maps on von Neumann algebras without central summands of type $I_1$.

\begin{corollary}\label{xxsec4.6}
Let $M$ be a von Neumann algebra without central summands of type $I_1$.
Then every $k$-commuting mapping of $\mathcal{M}$ is proper.
\end{corollary}

\subsection{Block upper and lower triangular matrix algebras}
\label{xxsec4.7}

Let $\mathbb{C}$ be the complex field. Let $\mathbb{N}$ be the set
of all positive integers and let $n\in \mathbb{N}$. For any
positive integer $m$ with $m\leq n$, we denote by $\bar{d}=(d_1,
\cdots, d_i, \cdots, d_m)\in \mathbb{N}^m$ an ordered $m$-vector
of positive integers such that $n=d_1+\cdots +d_i+\cdots+d_m$. The
\textit{block upper triangular matrix algebra}
$B^{\bar{d}}_n(\mathbb{C})$ is a subalgebra of $M_n(\mathbb{C})$
with form
$$
B^{\bar{d}}_n(\mathbb{C})=\left[
\begin{array}
[c]{ccccc}%
M_{d_1}(\mathbb{C})  & \cdots & M_{d_1\times d_i}(\mathbb{C}) & \cdots & M_{d_1\times d_m}(\mathbb{C})\\
& \ddots & \vdots &  & \vdots  \\
 & & M_{d_i}(\mathbb{C}) & \cdots & M_{d_i\times d_m}(\mathbb{C}) \\
 & O  &  & \ddots & \vdots \\
 &  &  & & M_{d_m}(\mathbb{C})  \\
\end{array}
\right].
$$
Likewise, the \textit{block lower triangular matrix algebra}
$B^{\prime \bar{d}}_n(\mathbb{C})$ is a subalgebra of
$M_n(\mathbb{C})$ with form
$$
B^{\prime \bar{d}}_n(\mathbb{C})=\left[
\begin{array}
[c]{ccccc}%
M_{d_1}(\mathbb{C})  &  &  &  & \\
\vdots& \ddots &  & O &   \\
M_{d_i\times d_1}(\mathbb{C}) & \cdots & M_{d_i}(\mathbb{C}) &  &  \\
 \vdots&   & \vdots & \ddots &  \\
M_{d_m\times d_1}(\mathbb{C}) & \ldots & M_{d_m\times d_i}(\mathbb{C}) & \ldots & M_{d_m}(\mathbb{C})  \\
\end{array}
\right]
$$
Note that the full matrix algebra $M_n(\mathbb{C})$ of all
$n\times n$ matrices over $\mathbb{C}$ and the upper(resp. lower)
triangular matrix algebra $T_n(\mathbb{C})$ of all $n\times n$
upper triangular matrices over $\mathbb{C}$ are two special cases
of block upper(resp. lower) triangular matrix algebras. If $n\geq
2$ and $B^{\bar{d}}_n(\mathbb{C})\neq M_n(\mathbb{C})$, then
$B^{\bar{d}}_n(\mathbb{C})$ is an upper triangular algebra and can
be written as
$$
B^{\bar{d}}_n(\mathbb{C})=\left[
\begin{array}
[c]{cc}%
B^{\bar{d}_1}_j(\mathbb{C}) & M_{j\times (n-j)}(\mathbb{C})\\
O_{(n-j)\times j} & B^{\bar{d}_2}_{n-j}(\mathbb{C})\\
\end{array}
\right],
$$
where $1\leq j < m$ and $\bar{d}_1\in \mathbb{N}^j, \bar{d}_2\in
\mathbb{N}^{m-j}$. Similarly, if $n\geq 2$ and
$B^{\prime\bar{d}}_n(\mathbb{C})\neq M_n(\mathbb{C})$, then
$B^{\prime \bar{d}}_n(\mathbb{C})$ is a lower triangular algebra
and can be represented as
$$
B^{\prime\bar{d}}_n(\mathbb{C})=\left[
\begin{array}
[c]{cc}%
B^{\prime\bar{d}_1}_j(\mathbb{C}) & O_{j\times (n-j)} \\
M_{(n-j)\times j}(\mathbb{C})& B^{\prime\bar{d}_2}_{n-j}(\mathbb{C})\\
\end{array}
\right],
$$
where $1\leq j < m$ and $\bar{d}_1\in \mathbb{N}^j, \bar{d}_2\in
\mathbb{N}^{m-j}$.

\begin{corollary}\label{xxsec4.7}
Every $k$-commuting mapping of the block upper triangular matrix
algebra $B^{\bar{d}}_n(\mathbb{C})$ {\rm (}resp. the block lower
triangular matrix algebra $B^{\prime \bar{d}}_n(\mathbb{C})${\rm
)} is proper.
\end{corollary}

\section{Topics for Future Potential Research}\label{xxsec5}

Although the main aim of this paper is to describe
the form of $k$-commuting mappings of  generalized matrix algebras,
the investigation of various additive mappings (associative-type, Jordan-type or Lie-type) on generalized
matrix algebras also have a great interest and should be further paid much attention.
The study of additive mappings on generalized matrix algebras is shedding light
on the investigation of functional identities in the background of such kind of algebras.
In the light of the motivation and contents of this article, we
will propose several topics with high potential and with merit for
future research in this field.

The theory of functional identities was initiated by Bre\v{s}ar at the beginning of 90's
in last century and it was greatly developed by Beidar, Bre\v{s}ar, Chebotar, and Martindale.
A functional identity (FI) of an algebra can be roughly described as an identical relation
involving arbitrary elements of the algebra together with functions. The scope of the
theory is to determine these functions or, in case this is not possible, to determine
the structure of the algebra admitting the given FI. The first functional identities
were introduced in the early 90's by Bre\v{s}ar as an attempt to unify several
results on centralizing mappings. Then the theory quickly developed through
a decade until reaching an ultimate stage that covers and unifies a number of
existing results. The main motivation for constructing a general theory
relies on the applications, and FI-theory has shown its strength in various areas. In
particular, it turned out to be the right tool in proving several conjectures formulated
by Herstein in 1961 concerning the description of Lie-type mappings in associative
rings \cite{Herstein}. For a full and nice account of the development of the theory of
functional identities and their applications, we refer the reader to the technical
literature \cite{BresarChebotarMartindale} for details.

Let $\mathcal{A}$ be an associative algebra. Let $F_1, F_2, G_1, G_2$ be mappings
from $\mathcal{A}$ into itself such that
$$
F_1(x)y+F_2(y)x+xG_2(y)+yG_1(x)=0 \eqno(5.1)
$$
for all $x, y\in \mathcal{A}$. This is a basic functional identity, which was one of the first functional
identities studied in prime algebras. The mappings $F_1, F_2, G_1$ and $G_2$ are looked on as unknowns
and the main purpose is to describe the form of these mappings. The functional identity $(5.1)$ is also
closely related to commuting mappings, which is due to the fact each
commuting additive mapping $F$ of $\mathcal{A}$ gives rise to the identity
$$
F(x)y+F(y)x-xF(y)-yF(x) = 0
$$
for all $x, y\in \mathcal{A}$. This identity is just one special case of $(5.1)$. Centralizing mappings
and commuting mappings can be considered as the most basic and important examples of
functional identities. However, the general theory of functional identities, which
was developed in \cite{BresarChebotarMartindale}, can not be applied to the context of
triangular rings, since these rings are not $d$-free. Nevertheless, Beidar, Bre\v{s}ar and Chebotar investigated certain functional
identities on upper triangular matrix algebras \cite{BeidarBresarChebotar}. Moreover, Cheung described
the form of commuting linear maps for a certain class of triangular algebras \cite{Cheung1, Cheung2}. Later, several
problems on certain types of mappings on triangular rings and algebras have been studied,
where some special examples of functional identities appear. Zhang and his students \cite{ZhangFengLiWu}
studied the functional identity of type (5.1) in the context of nest algebras. Han \cite{Han} considered the functional identity of type(5.1) on
CSL algebras and characterized the form of linear mappings $F_1, F_2, G_1, G_2\colon {\rm Alg}\mathcal{L}\longrightarrow \mathcal{M}$ satisfying $(5.1)$,
where $\mathcal{L}$ is a commutative subspace lattice generated by finite many
commuting independent nests on a complex separable Hilbert space $\mathbf{H}$
with ${\rm dim}\mathbf{H}\geq 3$, ${\rm Alg}\mathcal{L}$
is the CSL algebra associated with $\mathcal{L}$ and $\mathcal{M}$ be a $\sigma$-weakly
closed algebra containding ${\rm Alg}\mathcal{L}$. In two recent articles \cite{Eremita1, Eremita2},
Eremita studied functional identity $(5.1)$ in
triangular algebras. He succeeded in describing the form of additive mappings $F_1, F_2, G_1, G_2\colon \mathcal{T}\longrightarrow \mathcal{T}$
satisfying $(5.1)$ if a triangular ring $\mathcal{T}$ satisfies certain conditions. Moreover,
the notion of the
maximal left ring of quotients, which plays an important role in the study of functional identities
on (semi-)prime rings, is used to characterize those additive mappings $F_1, F_2, G_1, G_2$ \cite{Eremita2}.
It is predictable  that Eremita's approach, which is based on the notion of the maximal left ring
of quotients, enable us to generalize and unify a number of known results regarding mappings of
triangular algebras and generalized matrix algebras.

The functional identities in triangular algebras and full matrix algebras were
already studied in \cite{BeidarBresarChebotar, Eremita1, Eremita2}.
One would expect that the next step is to investigate
functional identities of generalized matrix algebras. The notion of generalized matrix algebras
efficiently unifies triangular algebras and full matrix algebras
together. The eventual goal of our systematic work is to deal
with all questions related to additive (or multiplicative) mappings of triangular
algebras and full matrix algebras under a unified frame, which is
the desirable generalized matrix algebras frame.

\begin{question}\label{xxsec5.1}
Let $\mathcal{G}=\mathcal{G}(A, M, N, B)$ be a generalized matrix algebra over a commutative ring $\mathcal{R}$ and let $F_1, F_2, G_1, G_2$ be mappings from $\mathcal{G}$ into itself such that
$$
F_1(x)y+F_2(y)x+xG_2(y)+yG_1(x)=0
$$
for all $x,y\in \mathcal{G}$. Describe the forms of $F_1, F_2, G_1, G_2$ satisfying the above condition.
\end{question}

Although people embark on studying functional identities in triangular algebras and full matrix algebras, the functional identities with additional structure has not been treated yet. For instance, the functional identities with automorphisms and derivations in generalized matrix algebras are worthy to be considered further.

\begin{proposition}{\rm(\cite[Proposition 4.2]{LiWei})}\label{xxsec5.2}
Let $\mathcal{G}=\mathcal{G}(A, M, N, B)$ be a generalized matrix algebra over a commutative ring $\mathcal{R}$. An $\mathcal{R}$-linear mapping $\Theta$ is a derivation of $\mathcal{G}$ if and only if
$\Theta$ has the form
$$
\Theta\left(\left[
\begin{array}
[c]{cc}%
a & m\\
n & b\\
\end{array}
\right]\right)
=\left[
\begin{array}
[c]{cc}%
\delta_1(a)-mn_0-m_0n & am_0-m_0b+\tau_2(m)\\
n_0a-bn_0+\nu_3(n) & n_0m+nm_0+\mu_4(b)\\
\end{array}
\right]
$$
where $m_0\in M, n_0\in N$ and
$$
\begin{aligned} \delta_1:& A \longrightarrow A, &
 \tau_2: & M\longrightarrow M, & \tau_3: & N\longrightarrow M,\\
\nu_2: & M\longrightarrow N, & \nu_3: & N\longrightarrow N , &
\mu_4: & B\longrightarrow B
\end{aligned}
$$
are all $\mathcal{R}$-linear mappings satisfying the following
conditions:
\begin{enumerate}
\item[{\rm(1)}] $\delta_1$ is a derivation of $A$ with
$\delta_1(mn)=\tau_2(m)n+m\nu_3(n);$

\item[{\rm(2)}] $\mu_4$ is a derivation of $B$ with
$\mu_4(nm)=n\tau_2(m)+\nu_3(n)m;$

\item[{\rm(3)}] $\tau_2(am)=a\tau_{2}(m)+\delta_1(a)m$ and
$\tau_2(mb)=\tau_2(m)b+m\mu_4(b);$

\item[{\rm(4)}] $\nu_3(na)=\nu_3(n)a+n\delta_1(a)$ and
$\nu_3(bn)=b\nu_3(n)+\mu_4(b)n.$
\end{enumerate}
\end{proposition}

In view of Proposition \ref{xxsec5.2} we have

\begin{question}\label{xxsec5.3}
Let $\mathcal{G}=\mathcal{G}(A, M, N, B)$ be a generalized matrix algebra over a commutative ring $\mathcal{R}$. Let $F_1, F_2, G_1, G_2$ be mappings from $\mathcal{G}$ into itself and $\Theta_1, \Theta_2, \Delta_1, \Delta_2$ be derivations from $\mathcal{G}$ into itself such that
$$
F_1(x)\Theta_2(y)+F_2(y)\Theta_1(x)+\Delta_1(x)G_2(y)+\Delta_2(y)G_1(x)=0
$$
for all $x,y\in \mathcal{G}$. What can we say about the forms of $F_1, F_2, G_1, G_2$ and those of $\Theta_1, \Theta_2, \Delta_1, \Delta_2$ satisfying the above condition ?
\end{question}

Let us next concentrate on the functional identities with automorphisms in generalized matrix algebras.
Automorphisms of generalized matrix algebras have been intensively considered
in \'Anh and van Wyk \cite{BobocDascalescuWyk}. Now we summarize some important
facts which are essential for our purposes.

\begin{proposition}{\rm(\cite[Theorem 3.6]{BobocDascalescuWyk})}\label{xxsec5.4}
Let $\mathcal{G}=\mathcal{G}(A, M, N, B)$ and $\mathcal{G}'=\mathcal{G}'(A', M',$ $ N', B')$ be two generalized matrix algebra over a commutative ring $\mathcal{R}$. Suppose that $A'$ and $B'$ have only trivial idempotents, and at least one of $M'$ and $N'$ is nonzero. Let $\Omega\colon \mathcal{G}\longrightarrow \mathcal{G}'$ be an $\mathcal{R}$-linear mapping. Then $\Omega$ is an isomorphism if and only if $\Omega$ has one of the following forms:
\begin{enumerate}
\item[(1)] $
\Omega \left(\left[ \begin{array}{cc} a & m\\
n & b \end{array}\right]\right) = \left[\begin{array}{cc}\gamma
(a) & \gamma (a) m^\prime_0 -m^\prime _0 \delta (b) + \mu(m) \\
n^\prime_0 \gamma (a)
- \delta (b) n^\prime_0 + \nu(n) & \delta (b) \end{array}\right]
$,
\end{enumerate}
where $\gamma\colon A\longrightarrow A^\prime$ and $\delta\colon B\longrightarrow B^\prime$ are two algebraic isomorphisms,
$\mu\colon M\longrightarrow M^\prime$ is a $(\gamma, \delta)$-bimodule isomorphism, \ $\nu\colon N\longrightarrow
N^\prime$ is a $(\delta, \gamma)$-bimodule isomorphism, $m^\prime_0 \in
M'$ and $n'_0 \in N'$ are two fixed elements.
\begin{enumerate}
\item[(2)] $\Omega \left(\left[ \begin{array}{cc} a & m \\
n & b \end{array}\right]\right) = \left[ \begin{array}{cc}\sigma
(b) & m'_\ast \rho(a) - \sigma (b) m'_\ast + \tau (n) \\
\rho(a)n'_\ast - n'_\ast \sigma (b) + \zeta (m) & \rho(a)
\end{array} \right]$,
\end{enumerate}
where $\rho\colon A
\longrightarrow B'$ and $\sigma\colon B \longrightarrow A'$ are two algebraic isomorphisms,
$\zeta\colon M\longrightarrow N'$ is a $(\rho, \sigma)$-bimodule isomorphism, $\tau\colon N\longrightarrow M'$
is a $(\sigma, \rho)$-bimodule isomorphism, $m'_\ast \in M'$ and $n'_\ast \in N'$ are two fixed elements.
\end{proposition}

A class of automorphisms established in the background of generalized matrix algebras permit us to
avoid the strong assumption---the diagonal algebras $A$ and $B$ in a generalized matrix
algebra $\mathcal{G}=\left[\smallmatrix A & M\\
N & B \endsmallmatrix \right]$ have only trivial idempotents---when we ideal with
those additive mappings with such kind of automorphisms on $\mathcal{G}$.
Moreover, their expression forms will be more simpler than the form of $\Omega$ in the above theorem.

Let $1$ (resp. $1^\prime$) be the identity of the algebra $A$ (resp.
$B$), and let $I$ be the identity of the generalized matrix algebra
$\mathcal{G}$. We will use the following notations:
$$
P=\left[
\begin{array}
[c]{cc}%
1 & 0\\
0 & 0\\
\end{array}
\right], \hspace{8pt} Q=I-P=\left[
\begin{array}
[c]{cc}%
0 & 0\\
0 & 1^\prime\\
\end{array}
\right]
$$
and
$$
\mathcal{G}_{11}=P{\mathcal G}P, \ \ \
\mathcal{G}_{12}=P{\mathcal G}Q, \ \ \
\mathcal{G}_{21}=Q\mathcal{G}P, \ \ \
\mathcal{G}_{22}=Q{\mathcal G}Q.
$$
Thus the generalized matrix algebra $\mathcal{G}$ can be written as
$$
\mathcal{G}=P{\mathcal G}P+P{\mathcal G}Q+Q\mathcal{G}P+Q{\mathcal G}Q
=\mathcal{G}_{11}+\mathcal{G}_{12}+\mathcal{G}_{21}+\mathcal{G}_{22}.
$$
Here, $\mathcal{G}_{11}$ and $\mathcal{G}_{22}$ are subalgebras of
$\mathcal{G}$ which are isomorphic to $A$ and $B$, respectively.
$\mathcal{G}_{12}$ is a $(\mathcal{G}_{11},
\mathcal{G}_{22})$-bimodule which is isomorphic to the $(A,
B)$-bimodule $M$. $\mathcal{G}_{21}$ is a $(\mathcal{G}_{22},
\mathcal{G}_{11})$-bimodule which is isomorphic to the $(B,
A)$-bimodule $N$. It should be remarked that $\pi_A(\mathcal{Z(G)})$
and $\pi_B(\mathcal{Z(G)})$ are isomorphic to $P\mathcal{Z(G)}P$ and
$Q\mathcal{Z(G)}Q$, respectively. Then there is an algebraic
isomorphism $\chi\colon P\mathcal{Z(G)}P\longrightarrow
Q\mathcal{Z(G)}Q$ such that $am=m\chi(a)$ for all $m\in
P\mathcal{G}Q$. There is also an algebraic
isomorphism $\varpi\colon Q\mathcal{Z(G)}Q\longrightarrow
P\mathcal{Z(G)}P$ such that $bn=n\varpi(b)$ for all $n\in
Q\mathcal{G}P$.

\begin{definition}\label{xxsec5.5}\cite[Definition 5.1.6]{Cheung1}
An automorphism $\xi$ of a generalized matrix algebra $\mathcal{G}=\left[\smallmatrix A & M\\
N & B \endsmallmatrix \right]$ is said to be \textit{partible} with respect to $A, B, M, N$, if it can be written as
$\xi=\phi_{r} \overline{\xi}$, where $r\in \mathcal{G}$, $\phi_{r}$ is an inner
automorphism $\phi_{r}(x)=r^{-1}xr$, which is induced by the element $r$, and
$\overline{\xi}$ is an automorphism of $\mathcal{G}$ satisfying the conditions
$\overline{\xi}(P\mathcal{G}P)=P\mathcal{G}P$, $\overline{\xi}(P\mathcal{G}Q)=P\mathcal{G}Q$,
$\overline{\xi}(Q\mathcal{G}P)=Q\mathcal{G}P$ and  $\overline{\xi}(Q\mathcal{G}Q)=Q\mathcal{G}Q$.

A generalized matrix algebra $\mathcal{G}=\left[\smallmatrix A & M\\
N & B \endsmallmatrix \right]$ is said to be $partible$
if every automorphism of $\mathcal{G}$ is partible.
\end{definition}

Indeed, there are a number of generalized matrix algebras which are partible, see \cite{FosnerLiangWei}.
Those triangular algebras what Mart\'{\i}n Gonz\'{a}lez, Repka and S\'anchez-Ortega considered in
\cite{MartinRepkaSanchez, RepkaSanchez, Sanchez} are exactly partible triangular algebras.
Of course, those algebras are also partible generalized matrix algebras.  They
take advantage of structural features of partible triangular algebras to describe commuting mapping
with automorphisms on these algebras. Comparing the old
version \cite{RepkaSanchez} and the new version \cite{MartinRepkaSanchez}, we observe that
partible triangular algebras not only cover upper triangular matrix algebras and Hilbert sapce nest
algebras, but also simplify the complicated computational process. Yu and Zhang \cite{YuZhang} investigated
commuting mappings with automorphisms on Hilbert space nest algebras. It should be remarked that
nest algebras on Hilbert spaces are partible, see \cite{FosnerLiangWei}.  In view of these works and
Proposition \ref{xxsec5.4} we ask

\begin{question}\label{xxsec5.6}
Let $\mathcal{G}=\mathcal{G}(A, M, N, B)$ be a partible generalized matrix algebra over a commutative ring $\mathcal{R}$.
Let $F_1, F_2, G_1, G_2$ be mappings from $\mathcal{G}$
into itself and $\Omega_1, \Omega_2, \Gamma_1, \Gamma_2$ be automorphisms from $\mathcal{G}$ into itself such that
$$
F_1(x)\Omega_2(y)+F_2(y)\Omega_1(x)+\Gamma_1(x)G_2(y)+\Gamma_2(y)G_1(x)=0
$$
for all $x,y\in \mathcal{G}$. How about are the forms of $F_1, F_2, G_1, G_2$ and those of $\Omega_1, \Omega_2, \Gamma_1, \Gamma_2$ satisfying the above condition ?
\end{question}

\vspace{4mm}
\noindent{\bf Acknowledgements}
\vspace{3mm}

\noindent Part of this research work was done when the third author visited the School of Mathematics
and Statistics at Beijing Institute of Technology in the winter of 2016. She takes this opportunity to
express his sincere thanks to the School of Mathematics and Statistics and the Office of International Affairs
at Beijing Institute of Technology for the hospitality extended to them during his visit.We are deeply grateful
to a special Training Program of International Exchange and Cooperation of the Beijing Institute of Technology.
We would like to thank the tremendous job of the anonymous referee, who, apart from a very thorough
report which helped to correct a number of minor errors, lacunae, and other inaccuracies (both mathematical
and pedagogical), also taught me some theory of generalized matrix algebras. And last but not least, we are
sincerely grateful to Professor Pham Ngoc ¨¢nh for his kind consideration and his advices on modifying our
manuscript.


\end{document}